\documentclass[12pt]{article}     
\usepackage[utf8]{inputenc}
\usepackage[russian]{babel}

\usepackage[T1]{fontenc}

\usepackage[warn]{mathtext}          
\usepackage[T2A]{fontenc}            

\newcommand*{\bfrac}[2]{\genfrac{}{}{0pt}{}{#1}{#2}}

\usepackage{graphicx}      
\usepackage{psfrag}        
\usepackage{caption2}      
\usepackage{soul}          
\usepackage{fancyhdr}      
\usepackage{multirow}      
\usepackage{ltxtable}      
\usepackage{paralist}      
\usepackage[perpage]{footmisc} 

\usepackage{amsmath}
\usepackage{amsfonts}
\usepackage{amssymb}
\usepackage{amsthm}
\usepackage{graphicx,color,framed} 

\DeclareMathOperator{\Res}{Re}
\DeclareMathOperator{\Ims}{Im}

\DeclareMathOperator{\mmod}{mod}
\newcommand{\Z}{\scriptstyle}
\newcommand{\D}{\displaystyle}

\renewcommand{\le}{\operatorname{\leqslant}}
\renewcommand{\ge}{\operatorname{\geqslant}}

\usepackage[a4paper, top=30mm, left=25mm, right=25mm, bottom=30mm]{geometry}

\newtheorem{lemma}{Lemma}
\newtheorem{theorem}{Theorem}
\newtheorem{remark}{Remark}

\newtheorem{corollary} {Corollary}

\newenvironment{myproof2} {\par\noindent{\bf Proof.}} {\hfill$\scriptstyle\blacksquare$}

\begin{document}

\begin{center}
\textsc{\textbf{\Large {On the mean values of some multiplicative functions on the short interval}}}
\footnote{This research was supported by the grant of Russian Fund of Fundamental Researches \No 12-01-31165.}
\end{center}

\begin{flushleft}
\begin{center}A.A. Sedunova\\ \end{center}
\end{flushleft}

{\bf Abstract.} In this paper we study the mean values of some multiplicative functions connected with the divisor function on the short interval of summation.
The asymptotic formulas for such mean values are proved.

\begin{section}{Introduction}
	
	In 1919, S.Ramanujan \cite{rmnj} announced the formula
	
	\begin{equation} \label{rmnj}
		\sum_{n \le X}{\frac{1}{\tau(n)}} = {\frac {X}{\sqrt{\ln X}} \left( A_0 + \frac{A_1}{\ln X} + \frac{A_2}{(\ln X)^2} + \ldots + \frac{A_N}{(\ln X)^N} + O \left( \frac{1}{(\ln X)^{N+1}}\right)\right)},
	\end{equation}
	where $A_j$ are some constants,
	$$A_0 = \frac{1}{\sqrt{\pi}} \prod_{p} { \sqrt{p (p-1)}\; \ln{\frac{p}{p-1}}},$$
	$\tau(n)$ denotes the number of divisors of $n$ and $N \ge 0$ is a fixed integer.	
	The complete proof of (\ref{rmnj}) was published in 1922 by B.M.Wilson \cite{wilson}.
	The general case (with $\tau_k(n)$ instead of $\tau(n)$) was considered by A.Ivi\'c \cite{ivic} in 1977.
	
	In this paper we generalize (\ref{rmnj}) and some other theorems of this type to the case when $n$ runs through the short interval of summation,
	i.e. the interval $x < n \le x+h$, where $x \rightarrow +\infty$ and $h \ll x^{\alpha}$ for a fixed $\alpha$, $0 < \alpha <1$.

	Suppose that $k \ge 2$ is fixed. The symbols $\sigma(n)$ and $r(n)$ stand for the sum of divisors of $n$ and for the number of representations of $n$
	by a sum of two squares $n = x^2 +y^2$, respectively.
	Let us define the multiplicative functions $f_j(n)$, $j=1,2,3,$ by the following relations:
	$$f_1(n) = \frac{1}{\tau_k(n)};\quad f_2(n) =  \D {\frac{\D \sigma(n)}{\D \tau(n)}};\quad f_3(n) =
			\left\{
			\begin{aligned}
				\frac{1}{r(n)},\;r(n) \neq 0;\\
				0,\;r(n) = 0.\\
			\end{aligned}
			\right.$$
	Finally, let
	$$S_j(x;h) = \sum_{x < n \le x+h}{f_j(n)}.$$
	Our goal is to prove the following theorems.
	
	\begin{theorem} \label{theorem1}
		Suppose that $N \ge 0$ is a fixed integer. Then the asymptotic formula
		$$S_1(x;h) = \sum_{x < n \le x+h}{\frac{1}{\tau_k(n)}} = \frac{h}{\left(\ln x\right)^{1 - \frac{1}{k}}} \left( A_0 + \frac{A_1}{\ln x} + \frac{A_2}{(\ln x)^2} + \ldots + \frac{A_N}{(\ln x)^N} + O \left(\frac {1}{(\ln x)^{N+1}} \right)\right),$$
		holds for $x \rightarrow +\infty$ and $h = x^{\alpha_k} e^{(\ln x)^{0.1}}$, where $\alpha_k=\frac{21k+5}{\Z 36k+5}$.
		
		Here the symbols $A_n = A_n(k)$ denote some positive constants that depend only on $n$ and $k$. In particular,
		$$A_0 = \left( \Gamma\left(\frac{1}{k}\right)\right) ^{-1} \prod_{p}{\left( 1 - \frac{1}{p}\right)^{\frac{1}{k}} F\left(1,1,k;\frac{1}{p}\right)},$$
		where
		$$F(a,b,c,z) = 1 + \frac{a\cdot b}{c}\,z + \frac{a(a+1)\cdot b(b+1)}{c(c+1)}\,z^{2}\,+\ldots\;$$
		is a hypergeometric function.
	\end{theorem}
	
	\begin{corollary}
	The formula
		$$\sum_{x < n \le x+h} {\frac{1}{\tau(n)}} = \frac{h}{\sqrt{\ln x}} \left( A_0 + \frac{A_1}{\ln x} + \frac{A_2}{(\ln x)^2} + \ldots + \frac{A_N}{(\ln x)^N} + O \left(\frac {1}{(\ln x)^{N+1}} \right)\right),$$
	holds true for any fixed $N \ge 0$ and for $h$ under the conditions
	$$x^{\frac{47}{77}}e^{(\ln x)^{0.1}} \le h \le x.$$
	\end{corollary}
	
	\begin{remark}
	The below table contains the approximate values of the constants $A_0(k)$ for $2 \le k \le 20$:
		\begin{center}
			\begin{tabular}{|l|c|l|c|}
				\hline
				{\it k} & $A_0(k)$ & {\it k} & $A_0(k)$\\ \hline
				2 & $0.546 855 96$ & 12 & $0.085 093 29$\\ \hline
				3 & $0.358 267 39$ & 13 & $0.078 426 06$\\ \hline
				4 & $0.264 796 54$ & 14 & $0.072 727 02$\\ \hline
				5 & $0.209 701 66$ & 15 & $0.067 799 67$\\ \hline
				6 & $0.173 497 45$ & 16 & $0.063 497 28$\\ \hline
				7 & $0.147 921 24$ & 17 & $0.059 708 09$\\ \hline
				8 & $0.128 902 38$ & 18 & $0.056 345 49$\\ \hline
				9 & $0.114 209 68$ & 19 & $0.053 341 30$\\ \hline
				10 & $0.102 519 63$ & 20 & $0.050 641 14$\\ \hline
				11 & $0.092 998 05$ & & \\
				\hline
			\end{tabular}
		\end{center}
	\end{remark}
	
	\begin{theorem} \label{theorem2}
		The formula
		\begin{equation} \label{f2}
		\sum_{x < n \le x+h}{\frac{\sigma(n)}{\tau(n)}} = \frac{h x}{\sqrt{\ln x}} \left( B_0 + \frac{B_1}{\ln x} + \frac{B_2}{(\ln x)^2} + \ldots + \frac{B_N}{(\ln x)^N} + O \left(\frac {1}{(\ln x)^{N+1}}\right)\right)
		\end{equation}
		holds true for any fixed $N \ge 0$ and for $h$ under the following conditions:
		$$x^{\frac{47}{77}} e^{(\ln x)^{0,1}} \le h \le x.$$
	\end{theorem}
	
	\begin{theorem} \label{theorem3}
		The formula
		\begin{equation} \label{f3}
		\sum_{x < n \le x+h}{\frac{1}{r(n)}} = \frac{h}{\left(\ln x\right)^{\frac{3}{4}}} \left( C_0 + \frac{C_1}{\ln x} + \frac{C_2}{(\ln x)^2} + \ldots + \frac{C_N}{(\ln x)^N} + O \left(\frac {1}{(\ln x)^{N+1}}\right)\right)
		\end{equation}
		holds true for any fixed $N \ge 0$ and for $h$ under the following conditions:
		$$x^{\frac{47}{77}} e^{(\ln x)^{0,1}} \le h \le x.$$
	\end{theorem}
	
	\begin{remark}
		The coefficients $B_n$ and $C_n$ depend only on $n$ with
		$$B_0 = \frac {1}{2\sqrt{\pi}} \prod_{p}{p \sqrt{\frac{p}{p-1}} \ln \left( 1+ \frac{1}{p}\right)} \approx 0.356903298$$
		and
		$$C_0 = \frac {2^{\frac{3}{4}}}{\Gamma\left(\tfrac{1}{4}\right)} \prod_{p \equiv 3 (\mmod 4)}{ \left( 1 - \tfrac{1}{p}\right)^{-\frac{3}{4}} \left( 1+\tfrac{1}{p}\right)^{-1}} \prod_{p \equiv 1 (\mmod 4)}{\left(1-\tfrac{1}{p}\right)^{\frac{1}{4}} \left( 1 + \tfrac{1}{2p} + \tfrac{1}{3p^2}+\ldots \right)} \approx $$
		$$ \approx	0.489330926.$$
	\end{remark}
	
	\begin{remark}
	The lower bounds for $h$ in the above theorems are not the best possible ones. They can be improved by using more precise upper bounds for $|\zeta(\sigma+it)|$
	in the strip $\frac{1}{2} \le \sigma \le 1$.
	\end{remark}
	
	{\bf Notations} \\
	In what follows, $C ,C_1, C_2\ldots$ denote positive absolute constants, which are, generally speaking, different in different relations.
	The symbol $(a,b)$ stands for the greatest common divisor of integer $a$ and $b$.
	Finally, $\theta, \theta_1, \theta_2, \ldots$ denote complex numbers with absolute values not greater than one, which are different in different relations.
\end{section}

\begin{section}{Auxilliary statements}
We need some auxilliary lemmas in order to prove theorems \ref{theorem1} - \ref{theorem3}.

\begin{lemma}
	Let $p$ be a prime number and let $\alpha \ge 1$. Then
	$$	\tau_k(p^{\alpha}) = C_{k+\alpha +1}^{k-1}, \;
		\sigma(p^{\alpha}) = \frac {p^{\alpha+1} - 1}{p -1}, \;
		r(p^\alpha) =
			\left\{
			\begin{aligned}
				4(\alpha + 1),\;if\;p=4k+1;\\
				4,\;if\;p=2;\\ 
				0,\;if\;p=4k+3\;and\;\alpha \;is\;an\;odd\;number;\\
				4,\;if\;p=4k+3\;and\;\alpha \;is\;an\;even\;number.\\
			\end{aligned}
			\right.
	$$
\end{lemma}
	
\begin{lemma}\label{perron}{\bf{(Perron's formula).}} 
	Suppose that the series $f(s) = \sum_{n=1}^{\infty}{a_n n^{-s}}$ converges absolutely for $\sigma > 1$, $|a_n| \le A(n)$, where
	$A(n)$ is a positive monotonicially increasing function of $n$ and
	$$\sum_{n=1}^{\infty}{|a_n| n^{-\sigma}} = O\left((\sigma-1)^{-\alpha}\right)$$
	for some $\alpha > 0$, as $\sigma \to  1+0$.
	Then the formula
	$$\sum_{n \le x}{a_n} = \frac{1}{2\pi i} \int_{b-iT}^{b+iT}{f(s) \frac{x^s}{s} ds} + O \left( \frac{x^b}{T(b-1)^{\alpha}}\right) + O \left( \frac{xA(2x)\ln x}{T}\right)$$
	holds true for any $b$, $1<b \le b_0$, $T \ge 2$, $x = N+\frac{1}{2}$ (the constants in O-symbols depend on $b_0$).
\end{lemma}
For a proof of the lemma, see \cite{vorkara}, pp. 334-336.

\begin{lemma}{\label{evzetasqr}}
	The estimate
	$$\int_{0}^{T}{\left|\zeta\left(\tfrac{1}{2} +it\right)\right|^2 dt} \ll T \ln T$$
	holds true for any $T \ge T_0 >1$.
\end{lemma}

\begin{myproof2}
	This lemma follows immediately from the theorem of Hardy and Littlewood (see, for example, \cite{titchmarch}, pp. 140-142).
\end{myproof2}

\begin{lemma} \label{evnofzeta}
	Let $\rho(u) = \frac{1}{2} - \{u\}$. Then the formula
	$$\zeta(s) = \frac{1}{2} + \frac{1}{s-1} +s\int_{1}^{\infty}{\frac{\rho(u) du}{u^{s+1}}},$$
	holds true for $s \neq 1$, $\Res s > 0$.
\end{lemma}
	For the proof, see \cite{vorkara}, pp. 24-25.

\begin{lemma} \label{T^()lnT}
	The estimates
	$$|\zeta(\sigma +it)| \ll t^{\frac{1 - \sigma}{3}} \ln t,\;\; \left| L(\sigma +it, \chi_4) \right| \ll t^{\frac{1 - \sigma}{3}} \ln t$$
	hold true for $|t| \ge t_0 >1$ and $\frac{\Z 1}{\Z 2} \le \sigma \le 1+ \frac{1}{\ln t}$.
\end{lemma}

\begin{myproof2}
This lemma can be easily derived from the approximate equations for $\zeta(s)$ and $L(s, \chi_4)$ (see, for example, \cite{lavrik} and \cite{titchmarch}, \S 7, Ch. IV) and from van der Corput's method of estimating of trigonometric sums.
\end{myproof2}
\vspace{0,2cm}

Let $N(\sigma, T)$ be the number of zeros of $\zeta(s)$ in the region $\Res s \ge \sigma, |\Ims s| \le T$.
Suppose that $q \ge 3$ is an integer and let $\chi$ be the Dirichlet's character modulo $q$.
Then the symbol $N(\sigma, T; \chi)$ stands for the number of zeros of the function $L(s, \chi)$ in the same domain.

\begin{lemma} \label{density}
	The estimates
	$$N(\sigma, T) \ll T^{\frac{12}{5} (1-\sigma)} (\ln T)^{44},$$
	$$\sum_{q \le  Q}\; \sideset{}{^*} {\sum}_{\chi \mmod \;Q} {N(\sigma,T; \chi)} \ll (Q^{2}T)^{\frac{12}{5}(1-\sigma)} (\ln QT)^{22}$$
	hold uniformly for $\frac{1}{\Z 2} \le \sigma \le 1$, $T \ge T_0$ and for $Q \ge 2$
	(the symbol $\sum ^*$ means the summation over all primitive characters $\chi$ modulo $q$).
\end{lemma}

\begin{lemma} \label{vinlim}
	There exist absolute positive constants $t_0$ and $C$ such that $\zeta(s) \neq 0$, $L(s, \chi_4) \neq 0$ in the region
	$$|t| \ge t_0, \;\; \sigma \ge 1 - \varrho(t),\; \varrho(t) = C (\ln \ln t)^{-\frac{1}{3}} (\ln t)^{-\frac{2}{3}}.$$
\end{lemma}

\end{section}

\begin{section}{Proof of the main results}
In this section we give the proofs of theorem 1, 2 and 3.
	\begin{subsection}{The mean-value of the function $\frac{\D 1}{\D \tau_k(n)}$ on the short interval}
		Suppose that $\sigma = \Res s > 1$ and let
		$$F(s) = \sum_{n=1}^{\infty}{\frac{1}{\tau_k(n)} \cdot n^{-s}}.$$
		This series converges absolutely, since 
		$$|F(s)| \le \sum_{n=1}^{\infty}{|a_n| \cdot n^{-\sigma}} \le \frac{1}{k} \sum_{n=1}^{\infty}{n^{-\sigma}} = 
		\frac{1}{k} \left( 1 + \int_{1}^{\infty}{\frac{du}{u^\sigma}}\right) = \frac{1}{k} \left( 1 + \frac{1}{\sigma-1}\right).$$
		
		Setting $a_n = \displaystyle{\frac{1}{\tau_k(n)}}, A(n) \equiv 1, b = 1 +\displaystyle{\frac{1}{\ln x}}, \alpha =1$ in lemma \ref{perron}, we get
		$$S_1=S(x,h;f_1)=I+O(R),$$
		where
		$$I=\frac{1}{2\pi i}\int_{b-iT}^{b+iT}{F(s) \frac{(x+h)^s-x^s}{s}ds},\;\;\; R=\frac{x^b}{T(b-1)}+\frac{xA(2x)\ln x}{T} \ll \frac{x \ln x}{T}.$$
		Further, $F(s) = \prod_p{F_p(s)}$, where
		$$F_p(s)=1 + \frac{1}{\tau_k(p)p^{s}} + \frac{1}{\tau_k(p^2)p^{2s}} + \ldots = 1+\frac{1!}{kp^s} + \frac{2!}{k(k+1)p^{2s}} + \ldots$$
		Writing $F_p(s)$ in the form
		$$F_p(s)=\left( 1 - \frac{1}{p^s}\right)^{-\frac{1}{k}} \left( 1 - \frac{1}{p^{2s}}\right)^{m_k} G_p(s),\quad{\rm where}\quad m_k=\frac{(k-1)^2}{2k^2(k+1)},$$
		we obtain
		$$F(s)=\frac{(\zeta(s))^{\frac{1}{k}}}{(\zeta(2s))^{m_k}} G(s),$$
		where
		$$G(s)=\prod_p{G_p(s)}=\prod_p{\left( 1 - \frac{1}{p^s}\right)^{-\frac{1}{k}} \left( 1 - \frac{1}{p^{2s}}\right)^{m_k} (1+u(s)+v(s))},$$
		$$u(s)=\frac{1}{kp^s}+\frac{2}{k(k+1)p^{2s}}, \;v(s)=\frac{3!}{k(k+1)(k+2)p^{3s}}+\frac{4!}{k(k+1)(k+2)(k+3)p^{4s}}+ \ldots$$
		Now we continue the function $F(s)$ to the left of the line $\Res s =1$.
		Suppose that $\frac{1}{2} \le \sigma \le 1$. Then the following estimates hold true:
		\begin{description}
			\item $$|u(s)| \le \frac{1}{kp^\sigma}\left( 1+\frac{2}{k+1}\frac{1}{p^\sigma}\right) \le \left(1 + \frac{2}{3\sqrt{2}}\right)\frac{1}{kp^\sigma} < \frac{3}{2} \frac{1}{kp^\sigma};$$
			\item $$|v(s)| \le \frac{1}{3kp^\sigma} \left( \frac{2\cdot3}{3\cdot4} + \frac{2\cdot3\cdot4}{3\cdot4\cdot5}\frac{1}{p^\sigma}+ \ldots \right) \le \frac{1}{2kp^{3\sigma}} \left( 1 + \frac{1}{p^\sigma} + \ldots \right) \le \frac{7}{4} \frac{1}{kp^\sigma};$$
			\item $$|u(s)+v(s)| \le \frac{3}{2k}\frac{1}{kp^\sigma}+\frac{7}{4k}\frac{1}{p^{3\sigma}} < \frac{5}{2}\frac{1}{kp^\sigma};$$
			\item $$|u(s) \cdot v(s)| \le \frac{3}{2k}\frac{7}{4k}\frac{1}{p^{4\sigma}} = \frac{21}{8k^2}\frac{1}{p^{4\sigma}} \le \frac{21}{16\sqrt{2}} \frac{1}{kp^{3\sigma}} < \frac{1}{kp^{3\sigma}}.$$
		\end{description}
		
		Now let us consider the expansion
		$$\ln(1+u(s)+v(s))=(u+v)-\frac{1}{2}(u+v)^2+\frac{1}{3}(u+v)^3 - \ldots$$
		Obviously, we have
		$$\left|\frac{1}{3}(u+v)^3 - \frac{1}{4}(u+v)^4+\ldots\right| \le \left(\frac{5}{2}\frac{1}{kp^\sigma}\right)^3+\frac{1}{4}\left(\frac{5}{2}\frac{1}{kp^\sigma}\right)^4 +\ldots \le $$
		$$\le \frac{1}{3}\left(\frac{5}{2}\frac{1}{kp^\sigma}\right)^3\frac{1}{1-\frac{5}{2kp^\sigma}} \le \frac{1}{3} \left( \frac{5}{2}\right)^3 \frac{1}{2^2} \left( 1 - \frac{5}{4\sqrt{2}}\right)^{-1} \frac{1}{kp^{3\sigma}} < \frac{45}{4} \frac{1}{kp^{3\sigma}}.$$
		Next,
		\begin{equation} \label{uv}
			(u+v)-\frac{1}{2}(u+v)^2=\left(u - \frac{u^2}{2}\right) + \left(v - \frac{v^2}{2}-uv \right).
		\end{equation}
		The second term on the right hand of (\ref{uv}) does not exceed in absolute value
		$$\frac{7}{4k} \frac{1}{p^{3\sigma}} + \frac{1}{2} \left( \frac{7}{4k}\right)^2\frac{1}{p^{6\sigma}} + \frac{1}{kp^{3\sigma}} < \frac{3}{kp^{3\sigma}}.$$
		Moreover,
		$$u-\frac{u^2}{2} = \frac{1}{kp^s}+\frac{2}{k(k+1)}\frac{1}{p^{2s}} -\frac{1}{2k^2p^{2s}} +\frac{3 \theta}{4}\frac{1}{kp^{3\sigma}} =
		\frac{1}{kp^s} + \frac{3k-1}{2k^2(k+1)}\frac{1}{p^{2s}} + \frac{3 \theta}{4} \frac{1}{kp^{3\sigma}}.$$
		Therefore,
		$$\ln (1+u(s)+v(s)) = \frac{1}{kp^s} + \frac{3k-1}{2k^2(k+1)}\frac{1}{p^{2s}} + \frac{15 \theta_1}{kp^{3s}}.$$
		Further,
		$$ \ln \left( 1 - \frac{1}{p^s}\right) = -\frac{1}{p^s}-\frac{1}{2p^{2s}}- \ldots = -\frac{1}{p^s}-\frac{1}{2p^{2s}} + \frac{7\theta_2}{6}\frac{1}{p^{3\sigma}},$$
		$$\ln \left( 1 - \frac{1}{2p^{2s}}\right) = -\frac{1}{p^{2s}} - \ldots = -\frac{1}{p^{2s}} + \frac{5\theta_3}{4} \frac{1}{p^{3\sigma}}.$$
		Hence
		$$\ln G_p(s) = \ln (1+u(s)+v(s)) +\frac{1}{k} \ln \left( 1 - \frac{1}{p^s}\right) - m_k \ln \left(1-\frac{1}{p^{2s}}\right) =$$
		$$ = \frac{1}{kp^s}+\frac{3k-1}{2k^2(k+1)}\frac{1}{p^{2s}} + \frac{3 \theta_1}{4} \frac{19}{kp^{3\sigma}} -\frac{1}{kp^s}-\frac{1}{2kp^{2s}} + \frac{7\theta_2}{6}\frac{1}{kp^{3\sigma}}
		+\frac{m_k}{p^{2s}} + \frac{5\theta_3}{4} \frac{m_k}{p^{3\sigma}}
		= \frac{16 \theta}{p^{3\sigma}}.$$
		
		Finally, for $\frac{1}{2} \le \sigma \le 1$ we get
		$$\left| \sum_{p}{\ln G_p(s)}\right| \le \frac{19}{k} \sum_{p}{\frac{1}{p^{\frac{3}{2}}}} < \frac{19}{k} < \frac{19}{2} = 9.5,$$
		$$-9.5 \le \ln |G(s)| \le 9.5, \;\;\; e^{-9.5} \le |G(s)| \le e^{9.5}.$$
		
		Let $\Gamma$ be the boundary of the rectangle with the vertices $\frac{\Z 1}{\Z 2} \pm iT, b \pm iT$, where the zeros of $\zeta(s)$ of the form
		$\frac{\Z 1}{\Z 2}+i\gamma$, $|\gamma|<T$, are avoided by the semicircles of the infinitely small radius lying to the right of the line $\Res s = \frac{\Z 1}{\Z 2}$,
		the pole of $\zeta(2s)$ at the point $s = \frac{\Z 1}{\Z 2}$ is avoided by two arcs $\Gamma_1$ and $\Gamma_2$ with the radius $\frac{\Z 1}{\Z \ln x}$,
		and let a horizontal cut be drawn from the critical line inside this rectangle to each zero $\rho = \beta + i\gamma$, $\beta > \frac{\Z 1}{\Z 2}$, $|\gamma| < T$.
		Then the function $F(s)$ is analytic inside $\Gamma$.
		By the Cauchy residue theorem,
		
		$$j_0 = - \sum_{k=1}^{8}{j_k} - \sum_{\rho}{j_{\rho}} = - (j_4+j_5) - \sum_{k \neq 4,5}{j_k} - \sum_{\rho}{j_{\rho}}.$$
		
		\begin{figure}[tbh]
			\begin{center}
				\includegraphics{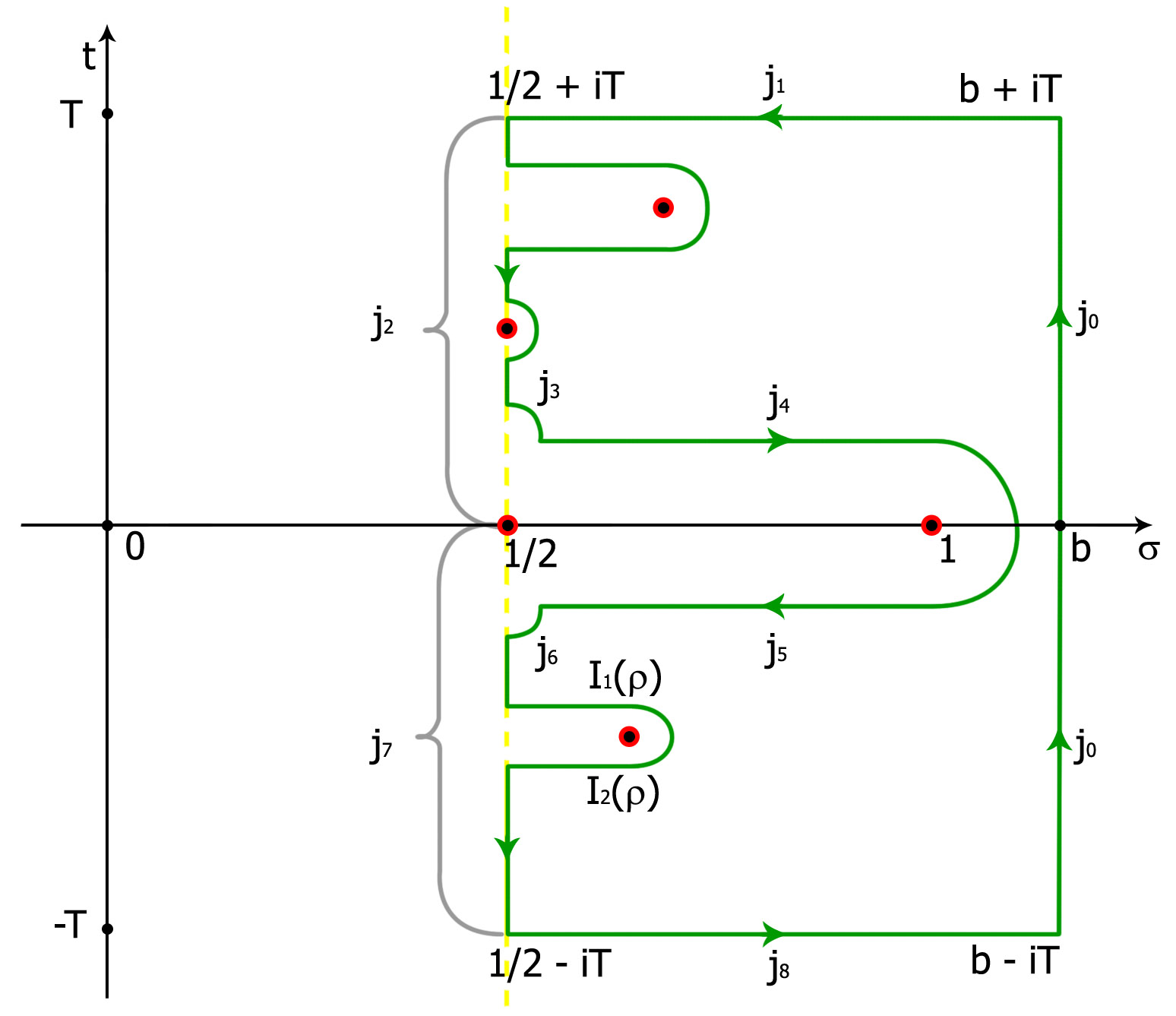}
			\end{center}
		\end{figure}
		
		By lemma \ref{T^()lnT},
		$$F(s) \ll T^{\frac{1-\sigma}{3k}} (\ln T)^{m_k} \ll T^{\frac{1-\sigma}{3k}} (\ln T)^{0.6}.$$
		Then
		$$|j_1|=\left| \frac{1}{2 \pi i} \int_{b + iT}^{\frac{1}{2} + iT}{F(s) \frac {(x + h)^s - x^s}{s} ds} \right| \ll
			\frac{1}{T} \int_{\frac{1}{2}}^{b}{T^{\frac{1 - \sigma}{3k}} \cdot (\ln x)^{0.6} x^{\sigma} d\sigma} \ll $$
		$$\ll \frac{x}{T} \int _{\frac{1}{2}}^{b}{\frac{x^{\sigma-1}}{T^{\frac{\sigma-1}{3k}}} (\ln x)^{0.6}} d\sigma \ll
		\frac{x}{T} \int _{\frac{1}{2}}^{b}{\left(\frac{x}{T^{\frac{1}{3k}}} \right)^{\sigma-1} (\ln x)^{0.6}} d\sigma \ll
			\frac{x}{T} (\ln x)^{0.6}.$$
		The similar estimate is valid for $j_8$.
		
		\begin{figure}[tbh]
			\begin{center}
				\includegraphics{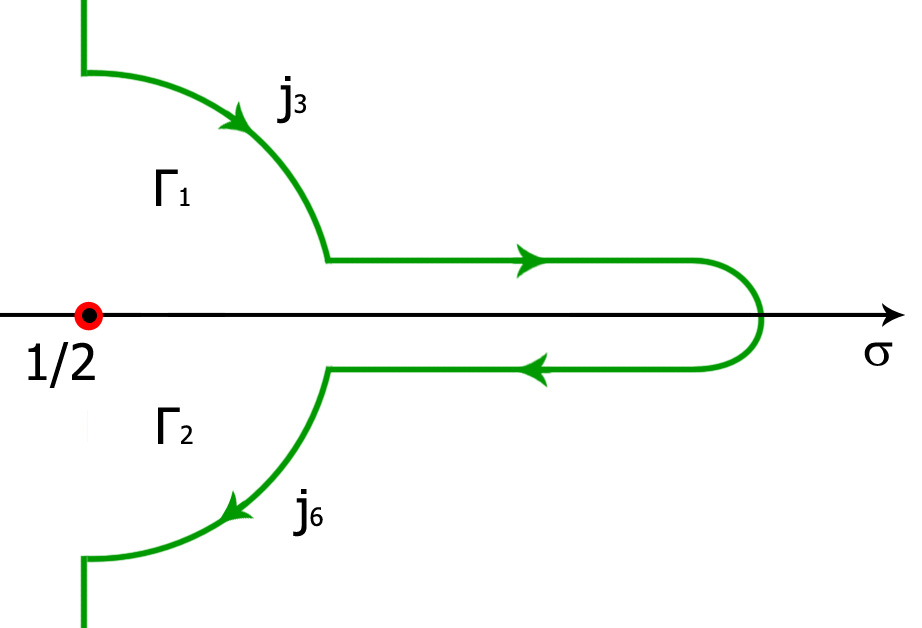}
			\end{center}
		\end{figure}
		
		By lemma \ref{evnofzeta}, on $\Gamma_1, \Gamma_2$ we have:
		$$|\zeta(s)| = \left| 0.5 + \frac {1}{s - 1} + s \int_{1}^{\infty}{\frac{\rho(u)}{u^{s + 1}} du} \right| \le 
			0.5 + 2 + 0.1 + 0.6 \cdot 0.5 \int_{1}^{\infty}{\frac {du}{u^{\frac{3}{2}}}} \le 3.2,$$
		$$\left| \zeta(2s)\right| \ge \frac {1}{|2s - 1|} - 1.1 \ge 0.5 \ln x - 1.1 > 0.4 \ln x.$$
		Hence
		$$|F(s)|  \le \left| G(s) \right| \frac {{3.2}^{\frac{1}{k}}}{(0.4 \cdot \ln x)^{m_k}} < 10.$$
		Therefore,
		$$|j_3 + j_6| \le \frac {1}{2\pi} \int_{\Gamma_1 \cup \Gamma_2}{|F(s)| \left| \frac {(x + h)^s - x^s}{s} \right| ds} \le
				\frac {10}{2\pi} \int_{- \frac{\pi}{2}}^{\frac{\pi}{2}}{\frac{2 \cdot (2x)^{\frac{1}{2} + \frac{1}{\ln x}}}{\frac{1}{2}} \cdot \frac {d\varphi}{\ln x}} \ll
				\frac{\sqrt{x}}{\ln x}.$$
		Further,
		$$|F(s)| \le |\zeta(s)|^{\frac{1}{k}}(\ln x)^{m_k}|G(s)| \ll (\ln x)^{m_k}\left|\zeta(\sigma+it)\right|^{\frac{1}{k}}.$$
		Hence
		$$|j_2| = \left| \text{p.v.} \frac{1}{2 \pi i}\int_{\frac{1}{2} + iT}^{\frac{1}{2} + \frac{i}{\ln x}}{F(s) \cdot \frac {(x + h)^s - x^s}{s} ds} \right|
			\ll \int_{\frac{1}{\ln x}}^{T} {(\ln x)^{m_k} \cdot \left|\zeta\left(\tfrac{1}{2} + it\right)\right|^{\frac{1}{k}}} \sqrt{x} \frac{dt}{t+1} \ll $$
			$$\ll (\ln x)^{m_k} \sqrt{x} \int_{0}^{T}{\left|\zeta\left(\tfrac{1}{2}+it\right)\right|^{\frac{1}{k}} \frac{dt}{t+1}} =
			(\ln x)^{m_k} \sqrt{x} \sum_{\nu \ge 0}{\int_{T/2^\nu}^{T/2^{\nu +1}}{\frac{|\zeta(\frac{1}{2} +it)|^{\frac{1}{k}}}{t+1}dt}}.$$
		Denoting the summands in the last sum by $j(\nu)$ and taking $X=T\cdot2^{-\nu}$, by the Hölder inequality we get:
		$$j(\nu) \ll \frac{1}{X} \left( \int_{X}^{2X}{\left|\zeta\left(\tfrac{1}{2}+it\right)\right|^2dt}\right)^{\frac{1}{2k}} X^{1-\frac{1}{2k}} \ll
		\frac{1}{X} \left( X\ln X\right)^{\frac{1}{2k}}X^{1-\frac{1}{2k}} \ll
		(\ln X)^{\frac{1}{2k}} \ll (\ln T)^{\frac{1}{2k}}.$$
		Hence,
		$$\sum_{\nu \ge 0}{j(\nu)} \ll (\ln T)^{1+\frac{1}{2k}} \ll (\ln T)^{\frac{5}{4}}.$$
		Then the upper bound for $j_2$ has the form
		$$|j_2| \ll (\ln x)^{m_k} \sqrt{x} (\ln x)^{\frac{5}{4}} = \sqrt{x} (\ln x)^{m_k+\frac{5}{4}}.$$
		The integral $j_7$ is estimated as above.
		
		The main term arises from the calculation of $j_4$ and $j_5$. Let us define the entire function $w(s)$ by the relation
		$$\zeta(s) = \frac{w(s)}{s-1}$$
		and let $s = 1 - u + i\cdot0$, where $0 \le u \le \frac{1}{2}$. Then
		$$\sqrt[k]{\zeta(s)} = \frac {\sqrt[k]{w(s)}}{\sqrt[k]{-u+i \cdot 0}}.$$
		Since $-u+i\varepsilon \to u\cdot e^{\pi i}$ as $\varepsilon \to +0$, then
		$$\sqrt[k]{-u+i \cdot 0} = \sqrt[k]{u} e^{\frac{\pi i}{k}}, \;\; \sqrt[k]{\zeta(s)}=\frac{\sqrt[k]{w(\sigma)}}{\sqrt[k]{u}} e^{-\frac{\pi i}{k}}.$$
		Therefore, on the upper edge of the cut we have
		$$F(s) = \frac{\sqrt[k]{w(1-u)}}{(\zeta(2-2u))^{m_k}}\; G(1-u)\frac{e^{-\frac{\pi i}{k}}}{\sqrt[k]{u}} = \frac{\Pi(u)e^{-\frac{\pi i}{k}}}{\sqrt[k]{u}},$$
		where
		$$\Pi(u) = G(1-u) \frac{\sqrt[k]{w(1-u)}}{(\zeta(2-2u))^{m_k}}.$$
		
		Hence,
		$$j_4 = \frac{1}{2 \pi i} \int_{\frac{1}{2} + \frac{1}{\ln x} + i \cdot 0}^{1+i \cdot 0}{F(\sigma + i \cdot 0) \frac{(x+h)^s- x^s}{s} ds} = $$
		$$=\frac{1}{2\pi i} \int_{\frac{1}{2}+\frac{1}{\ln x}+i \cdot 0}^{1+i \cdot 0}{F(\sigma+i \cdot 0) \int_{0}^{h}{(x+u)^{s-1}du}\;ds} =$$
		$$= \frac{1}{2\pi i}\int_{x}^{x+h}{\int_{\frac{1}{2}+\frac{1}{\ln x}}^{1}{F(\sigma+i \cdot 0)y^{\sigma-1}d\sigma}\;dy}=
		\frac{e^{-\frac{\pi i}{k}}}{2\pi i} \int_{x}^{x+h}{{\int_{0}^{\frac{1}{2}-\frac{1}{\ln x}}{\frac{\Pi(u)y^{-u}}{\sqrt[k]{u}}}du}\;dy}.$$
		Suppose that $N \ge 0$ is fixed. Then
		$$\Pi(u)= \Pi_0+\Pi_1u + \Pi_2u^2 + \ldots + \Pi_N u^N + O_N(u^{N+1}),$$
		where
		$$\Pi_0=\Pi(0)=\frac{\sqrt[k]{w(1)}}{(\zeta(2))^{m_k}} \; G(1) = \frac{\prod_{p} \left(1-\frac{1}{p}\right)^{\frac{1}{k}} \left( 1-\frac{1}{p^2}\right)^{m_k} \left( 1 + \frac{1}{kp} + \ldots \right)}{(\zeta(2))^{m_k}} = $$
		$$=\prod_p{\left( 1 - \frac{1}{p}\right)^{\frac{1}{k}}\left( 1+ \frac{1}{kp} + \ldots\right)}.$$
		Thus, we have
		$$j_4 = \frac{e^{-\frac{\pi i}{k}}}{2\pi i} \int_{x}^{x+h}{\left(\sum_{0 \le n \le N}{\Pi_n \int_{0}^{\frac{1}{2}-\frac{1}{\ln x}}{\frac{u^{n}y^{-u}}{\sqrt[k]{u}}du}}+ O\left( J\right)\right)dy},$$
		where
		$$J = \int_{0}^{\frac{1}{2}-\frac{1}{\ln x}}{\frac{u^{N+1}y^{-u}}{\sqrt[k]{u}}du} \le  \frac{\Gamma\left(N+2-\frac{1}{k}\right)}{(\ln y)^{N+2-\frac{1}{k}}}.$$
		Using the estimate
		$$\int_{\lambda}^{\infty}{w^{k-\gamma}e^{-w}dw < ek!\;\lambda^{k-\gamma}e^{-\lambda}},$$
		where $\lambda>1$, $0 < \gamma <1$, $k \ge 1$, we easily get
		$$\int_{0}^{\frac{1}{2}-\frac{1}{\ln x}}{\frac{u^n y^{-u}}{\sqrt[k]{u}}du} =\frac{1}{(\ln y)^{n+1-\frac{1}{k}}} \int_{0}^{\ln y \left(\frac{1}{2}-\frac{1}{\ln x} \right)}{e^{-w} w^{n-\frac{1}{k}} dw}=$$
		$$=\frac{1}{(\ln y)^{n+1-\frac{1}{k}}} \left( \int_{0}^{\infty}{e^{-w} w^{n-\frac{1}{k}} dw} + \frac{\theta e n! (\ln y)^{n-\frac{1}{k}}}{\sqrt{y}}\right)=
		\frac{\Gamma\left( n+1-\frac{1}{k}\right)}{(\ln y)^{n+1-\frac{1}{k}}}+ \frac{\theta e n!}{\sqrt{y} \ln y}.$$
		Therefore,
		$$j_4 = \frac{e^{-\frac{\pi i}{k}}}{2\pi i} \int_{0}^{h}{\sum_{0 \le n \le N}{\frac{\Pi_n \Gamma(n+1-\frac{1}{k})}{(\ln (x+u))^{n+1-\frac{1}{k}}}}du}+O\left( \frac{h}{\left(\ln x\right)^{N+2-\frac{1}{k}}}\right).$$
		Let
		$$\varphi(x) = \frac{1}{\left( \ln x\right)^{n+1-\frac{1}{k}}}.$$
		Then the Lagrange mean-value theorem yields
		$$\varphi(x+u) = \varphi(x) + u \varphi^{'}(x+\theta_1 u) = \frac{1}{\left( \ln x\right)^{n+1-\frac{1}{k}}} + \frac{\theta_2 h \left( n+1-\frac{1}{k}\right)}{x \left( \ln x\right)^{n+2-\frac{1}{k}}}.$$
		Thus we get
		$$j_4 = \frac{h e^{-\frac{\pi i}{k}}}{2\pi i} \sum_{0 \le n \le N}{\frac{\Pi_n \Gamma\left( n+1-\frac{1}{k}\right)}{(\ln x)^{n+1-\frac{1}{k}}}}+ O\left( \frac{h}{(\ln x)^{N+2-\frac{1}{k}}}\right) + O\left( \frac{h^2}{x} \frac{1}{(\ln x)^{N+2-\frac{1}{k}}}\right),$$
		$$j_5 = -\frac{h e^{\frac{\pi i}{k}}}{2\pi i} \sum_{0 \le n \le N}{\frac{\Pi_n \Gamma\left( n+1-\frac{1}{k}\right)}{(\ln x)^{n+1-\frac{1}{k}}}}+ O\left( \frac{h}{(\ln x)^{N+2-\frac{1}{k}}}\right) + O\left( \frac{h^2}{x} \frac{1}{(\ln x)^{N+2-\frac{1}{k}}}\right).$$
		
		Finally,
		$$-(j_4+j_5)=
		-\frac{h}{(\ln x)^{1-\frac{1}{k}}} \left( \sum_{0 \le n \le N}{\frac{(-1)^n \Pi_n}{\Gamma \left( \frac{1}{k}-n\right) (\ln x)^n}} + O\left( \frac{1}{(\ln x)^{N+1}}\right) + O \left( \frac{h}{x(\ln x)^{N+1}}\right)\right).$$
		
		It remains to estimate the sum
		$$\sum_{|\gamma|<T}{j_{\rho}},\;j_{\rho}=I_1(\rho)+I_2(\rho),\;\text{where}\;\;\rho = \beta +i\gamma.$$
		
		\begin{figure}[tbh]
			\begin{center}
				\includegraphics{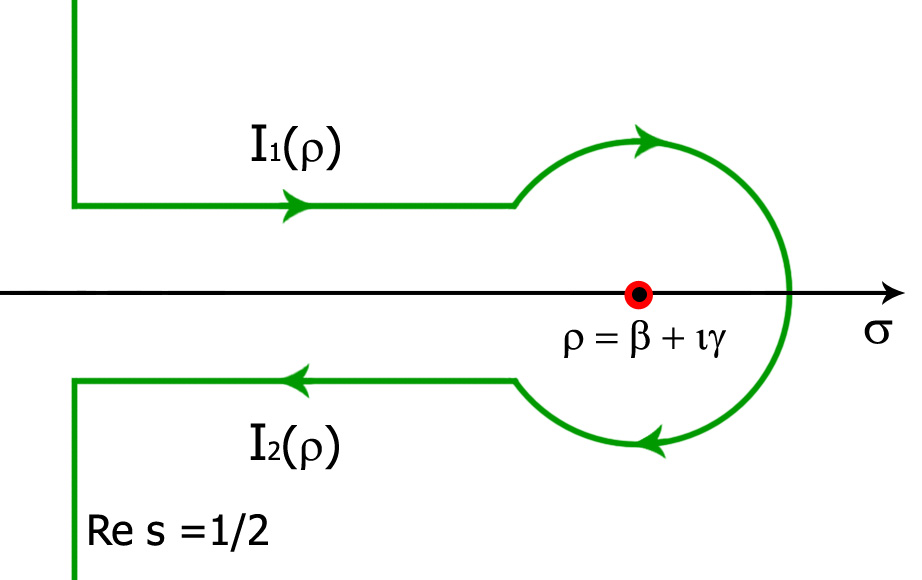}
			\end{center}
		\end{figure}
		Since
		$$\left| \frac{(x+h)^s-x^s}{s}\right| = \left| \int_{0}^{h}{(x+u)^{s-1}du}\right| \ll  \int_{0}^{h}{(x+u)^{\sigma-1}du} \ll h x^{\sigma-1},$$
		then
		$$I_1(\rho) \ll \int_{\frac{1}{2}}^{\beta}{(\ln x)^{m_k} |\zeta(\sigma+i\gamma)|^{\frac{1}{k}} h x^{\sigma-1} d\sigma} \ll
		\frac{h}{x} (\ln x)^{m_k} \int_{\frac{1}{2}}^{\beta}{x^{\sigma}  |\zeta(\sigma+i\gamma)|^{\frac{1}{k}} d\sigma}$$
		and the same estimate is valid for $I_2(\rho)$.
		Hence,
		$$|j_\rho| \ll \int_{\frac{1}{2}}^{\beta}{h x^{\sigma-1} (\ln x)^{m_k} T^{\frac{1-\sigma}{3k}} (\ln x)^{\frac{1}{k}} d\sigma} \ll
		h (\ln x)^{m_k+\frac{1}{k}} \int_{\frac{1}{2}}^{\beta}{\left( \frac{T^{\frac{1}{3k}}}{x}\right)^{1-\sigma} d\sigma} \ll$$
		$$\ll h (\ln x)^{m_k+\frac{1}{k}} \int_{\frac{1}{2}}^{1}{g(\rho,\sigma) \left( \frac{T^{\frac{1}{3k}}}{x}\right)^{1-\sigma} d\sigma},$$
		where
		$$g(\rho,\sigma) = \left\{
			\begin{aligned}
				1,\;\text{if}\;\sigma \le \beta,\\ 
				0,\;\text{if}\;\sigma > \beta.\\
			\end{aligned}
			\right.$$
		Applying lemma \ref{density}, we get
		$$\sum_{|\gamma|<T}{j_\rho} \ll h (\ln x)^{m_k+\frac{1}{k}} \int_{\frac{1}{2}}^{1-\varrho (T)}{\left( \sum_{|\gamma| < T}{g(\rho; \gamma)} \right) \left( \frac{T^{\frac{1}{3k}}}{x}\right)^{1-\sigma} d\sigma} \ll$$
		$$\ll h (\ln x)^{m_k+\frac{1}{k}} \int_{\frac{1}{2}}^{1-\varrho (T)} {N(\sigma; T) \left( \frac{T^{\frac{1}{3k}}}{x}\right)^{1-\sigma} d\sigma}\ll$$
		$$\ll h (\ln x)^{m_k+44+\frac{1}{k}} \int_{\frac{1}{2}}^{1-\varrho (T)} {\left( \frac{T^{\frac{1}{3k}}}{x}\right)^{1-\sigma} T^{\frac{12}{5}(1-\sigma)} d\sigma} \ll$$
		$$\ll  h (\ln x)^{45} \int_{\frac{1}{2}}^{1-\varrho (T)} {\left( \frac{T^{\frac{12}{5}+\frac{1}{3k}}}{x}\right)^{1-\sigma} d\sigma} \ll
		 h (\ln x)^{45} \left( \frac{T^{\frac{12}{5} + \frac{1}{3k}}}{x}\right)^{1-\varrho (T)}.$$
		Let $D(x)=e^{C_1(\ln x)^{0.8}}.$
		Choosing $T$ from the equation,
		$$T^{\frac{36k+5}{15k}}=x D^{-1}(x)$$
		we easily conclude that
		$$T=x^{\frac{15k}{36k+5}} D(x)^{-\frac{15k}{36k+5}}$$
		and
		$$j_0 = \frac{h}{(\ln x)^{1-\frac{1}{k}}} \left( \sum_{0 \le n \le N}{\frac{(-1)^n \Pi_n}{\Gamma(\frac{1}{k}-n) (\ln x)^n}}\right) + O(J),$$
		where
		
		\begin{multline}
			J = \frac{x}{T} (\ln x)^{0.6} + \frac{\sqrt{x}}{\ln x} + \sqrt{x} (\ln x)^{m_k + \frac{5}{4}} + \frac{h}{(\ln x)^{N+2-\frac{1}{k}}} + \frac{h^2}{x(\ln x)^{N+2-\frac{1}{k}}} + \\ + h(\ln x)^{45}  \left(\frac{x^{\frac{1}{3k}}}{D^{1+\frac{1}{3k}}}\right)^{1-\varrho (T)} \ll \\
			\ll \frac{x}{T} (\ln T)^{\frac{17}{20}} + \frac{h}{(\ln x)^{N+2-\frac{1}{k}}} + \frac{h^2}{x(\ln x)^{N+2-\frac{1}{k}}} + h(\ln x)^{45}  \left( \frac{x^{\frac{1}{3k}}}{D^{1+\frac{1}{3k}}}\right)^{1-\varrho (T)}.
		\end{multline}
		Obviously, the formula for $S_1$ is asymptotic if
		$$h \gg \frac{x}{T} (\ln x)^2 = x^{\alpha_k} e^{C_2 (\ln x)^{0.8}},$$	
		where
		$$\alpha_k = 1-\frac{15k}{36k+5}=\frac{21k+5}{36k+5}.$$
	\end{subsection}
	
	\begin{subsection}{The mean-value of the function $\frac{\D \sigma(n)}{\D \tau(n)}$ on the short interval}
	
		Suppose that $\sigma = \Res s > 2$ and let
		$$F(s) = \sum_{n=1}^{\infty}{\frac{\sigma(n)}{\tau_k(n)} \cdot n^{-s}}.$$
		This series converges absolutely, since 
		$$|F(s)| \le \sum_{n=1}^{\infty}{n^{1-\sigma}} \le  1 + \int_{1}^{\infty}{u^{1-\sigma}du} = 1+ \frac{1}{\sigma-2 }< +\infty $$
		Setting
		$$a_n = \displaystyle{\frac{\sigma(n)}{\tau_k(n)}},\;\;n \equiv A(n), b = 2 +\displaystyle{\frac{1}{\ln x}} \le \displaystyle{\frac{21}{10}},\;\;\alpha = 1$$
		in lemma \ref{perron}, we get:
		$$a_n= \frac{1}{\tau(n)} \sum_{d|n}{d} \le \frac{n \tau(n)}{\tau(n)} \le A(n),$$
		$$S_2=S(x,h;f_2)=I+O(R),$$
		where
		$$I=\frac{1}{2\pi i}\int_{b-iT}^{b+iT}{F(s) \frac{(x+h)^s- x^s}{s}ds},\;R=\frac{x^b}{T(b-1)}+\frac{xA(2x)\ln x}{T} \ll \frac{x^2 \ln x}{T}.$$
		Further,
		$$F(s) = \prod_p{F_p(s)},$$
		where
		$$F_p(s)=1 + \sum_{k=1}^{\infty}{\frac{1}{k+1} \cdot \frac{1}{p^{k(s-1)}} \cdot \frac{1-\frac{1}{p^{k+1}}}{1-\frac{1}{p}}} =
		1 + \frac{1}{2p^s} + \frac{1}{2p^{s-1}} + \sum_{k=2}^{\infty}{\frac{1}{k+1} \cdot \frac{1}{p^{k(s-1)}} \cdot \frac{1-\frac{1}{p^{k+1}}}{1-\frac{1}{p}}}.$$
		Obviously,
		$$F_p(s+1)=1+\frac{p+1}{2p^{s+1}}+\frac{p^2+p+1}{3p^{2(s+1)}} + \sum_{k=3}^{\infty}{\frac{p^{k+1}}{p-1}} \cdot \frac{p^{-(s+1)k}}{k+1} =
		1+u(s) + v(s),$$
		where
		$$v(s) = \frac{1}{2p^s}+\frac{1}{3p^{2s}},$$
		$$u(s)=\frac{1}{2p^{s+1}} +\frac{p+1}{3p^{2(s+1)}} + \sum_{k=3}^{\infty}{\frac{p^{k+1}}{p-1}} \cdot \frac{p^{-(s+1)k}}{k+1}.$$
		Suppose that $p \ge p_0=5$. Then
		$$|u(s)| \le \frac{1}{2p^{\sigma+1}} + \frac{p+1}{3p^{2(\sigma+1)}} + \sum_{k=3}^{\infty}{\frac{(k+1)p^k}{k+1} \frac{1}{p^{k(\sigma+1)}}} =
		\frac{1}{2p^{\sigma+1}} + \frac{p+1}{3p^{2(\sigma+1)}} + \sum_{k=3}^{\infty}{\frac{1}{p^{k\sigma}}} =$$
		$$= \frac{1}{p^{\sigma+1}} \left( \frac{1}{2} + \frac{1}{3p^\sigma}+\frac{1}{3p^{\sigma+1}} + \frac{1}{p^{2\sigma-1}} \cdot \frac{1}{1-\frac{1}{p^\sigma}}\right) \le
		\frac{C_1}{p^{\sigma+1}},$$
		where
		$$C_1=\frac{1}{2} + \frac{1}{3\sqrt{p_0}} \left( 1+ \frac{1}{p_0}\right) + \frac{1}{1-\frac{1}{\sqrt{p_0}}},$$
		
		$$|v(s)| \le \frac{1}{2p^\sigma}+\frac{1}{3p^{2\sigma}} = \frac{1}{p^\sigma} \left( \frac{1}{2} + \frac{1}{3p^\sigma}\right) \le \frac{C_2}{p^\sigma},\;C_2=\frac{1}{2}+\frac{1}{3\sqrt{p_0}}.$$
		Therefore,
		$$|u(s)+v(s)| \le \frac{C_1}{p^{\sigma+1}} + \frac{C_2}{p^{\sigma}} = \frac{C_0}{p^\sigma},$$
		$$C_0 = \frac{C_1}{p_0}+ C_2=\frac{1}{2}\left(1+\frac{1}{p_0}\right)+\frac{1}{3\sqrt{p_0}} \left(1+\frac{1}{p_0}+\frac{1}{p_0^2}\right)+\frac{1}{p_0-\sqrt{p_0}},$$
		$$|u(s)\cdot v(s)| \le \frac{C_1\cdot C_2}{p^{2\sigma+1}},\;\;\;|u(s)|^2 < \frac{C_1^2}{p^{2\sigma+2}},$$
		where
		$$C_0=\frac{C_1}{p_0}+C_2 = 1.1466...,\; C_1 = 2.4879...,\; C_2 = 0.6490...$$
		Furthermore,
		$$F_p(s+1)=\left( (v+u)-\frac{1}{2}(u^2+2uv+v^2) \right) +\sum_{k \ge 3}{\frac{(-1)^k}{k} (u+v)^k} = U+V,$$
		$$|V|=\left| \sum_{k \ge 3}{\frac{(-1)^k}{k} (u+v)^k}\right| \le \sum_{k \ge 3}{\frac{1}{k} \left( \frac{C_0}{p^\sigma}\right)^k} \le
		\frac{1}{3} \left(\frac{C_0}{p^\sigma}\right)^3 \sum_{k=0}^{\infty}{\left(\frac{C_0}{\sqrt{p_0}}\right)^k} = \frac{C_3}{p^{3\sigma}},$$
		where
		$$C_3=\frac{C_0^3}{3}\cdot \frac{1}{1-\frac{C_0}{\sqrt{p_0}}}=1.0314...$$
		Representing $U$ in the form
		$$U = (u+v)-\frac{u^2}{2}-uv-\frac{v^2}{2}=\frac{1}{2^s}+\frac{1}{3p^{2s}}+u-\frac{u^2}{2}-uv-\frac{1}{2}\left(\frac{1}{4p^{2s}}+\frac{1}{3p^{3s}}+\frac{1}{9p^{4s}}\right)=$$
		$$=\frac{1}{2p^s}+\frac{5}{24}\cdot \frac{1}{p^{2s}}+W,$$
		we have:
		$$|W| = \left| u - \frac{u^2}{2} -uv -\frac{1}{2} \left( \frac{1}{3p^{3s}}+\frac{1}{9p^{4s}}\right)\right| \le \frac{C_1}{p^{\sigma+1}} + \frac{C_1^2}{2} \cdot \frac{1}{p^{2\sigma+2}} +\frac{C_1 C_2}{p^{2\sigma+1}}+\frac{1}{2} \left( \frac{1}{3p^{3\sigma}} + \frac{1}{9p^{4\sigma}}\right) \le $$
		$$\le \frac{1}{p^{\sigma+1}} \left( C_1 + \frac{C_1^2}{2}\cdot \frac{1}{p_0^{\sigma+1}}+\frac{C_1 C_2}{p_0^{2\sigma+1}} +\frac{1}{2} \left( \frac{1}{3p_0^{2\sigma-1}} + \frac{1}{9p_0^{3\sigma-1}} \right) \right) \le $$
		$$\frac{1}{p^{\sigma+2}} \left( C_1 + \frac{C_1^2}{2} \cdot \frac{1}{p_0 \sqrt{p_0}} + \frac{C_1 C_2}{\sqrt{p_0}}+ \frac{1}{2} \left( \frac{1}{3} + \frac{1}{9\sqrt{p_0}}\right)\right) = \frac{C_4}{p^{\sigma+1}},$$
		where
		$$C_4 = C_1 + \frac{1}{6} + \frac{1}{\sqrt{p_0}} \left( \frac{C_1^2}{2p_0}+C_1 C_2 \frac{1}{18}\right) = 1.7297...$$
		Thus, for $p \ge p_0$ we obtain
		$$ \ln F_p(s+1) = \ln (1+v(s)+u(s)) = \frac{1}{2p^s} + \frac{5}{24p^{2s}} + \theta_1 \left( \frac{C_3}{p^{3\sigma}}+\frac{C_4}{p^{\sigma+1}}\right) = $$
		$$=\frac{1}{2p^s}+\frac{5}{24p^{2s}} + \frac{\theta_2 C_5}{p^{\sigma+1}},\;\;\;C_5=C_3+C_4=2.7612...$$
		Now let
		$$G_p(s) = F_p(s) \left( 1 - \frac{1}{p^{s-1}}\right)^{\frac{1}{2}} \left( 1-\frac{1}{p^{2(s-1)}}\right)^{-\frac{1}{24}}.$$
		Since
		$$\frac{1}{2} \ln \left(1-\frac{1}{p^s}\right) = -\frac{1}{2p^s}-\frac{1}{4p^{2s}} + \frac{\theta_3 C_6}{p^{3s}},\; C_6 = \frac{1}{6} \cdot \frac{1}{1-\frac{1}{\sqrt{p_0}}} = \frac{5+\sqrt{5}}{24},$$
		$$-\frac{1}{24} \ln \left( 1 - \frac{1}{p^{2s}}\right) = \frac{1}{24p^{2s}} + \frac{5}{8 \cdot 24} \cdot \frac{\theta_4}{p^{4\sigma}},$$
		then we finally get:
		
		$$|\ln G_p(s+1)| \le \frac{C_5}{p^{\sigma+1}} + \frac{C_6}{p^{3\sigma}} + \frac{5}{8 \cdot 24} \cdot \frac{1}{p^{4\sigma}} = 
		\frac{1}{p^{\sigma+1}} \left( C_5 + \frac{C_6}{p^{2\sigma-1}} + \frac{5}{8 \cdot 24} \cdot \frac{1}{p^{3\sigma-1}}\right) \le \frac{C_7}{p^{\sigma+1}},$$
		$$C_7 = C_5 + C_6 + \frac{5}{8 \cdot 24} \cdot \frac{1}{\sqrt{p_0}} = 3.2387...$$
		Thus,
		$$|\ln G_p(s+1)| \le \frac{C_7}{p^{\sigma+1}},\;\;\;C_7 = 3.2387...$$
		
		It remains to check the cases of $p=2$ and $p=3$.
		For $p=2$ we have
		$$F_2(s+1)=\sum_{k=0}^{\infty}{\frac{2^{k+1}-1}{k+1} {2^{-k(s+1)}}} = 2^{s+1} \sum_{k=0}^{\infty}{\frac{2^{k+1}-1}{k+1} {2^{-(k+1)(s+1)}}} =$$
		$$=2^{s+1}\left( \sum_{k=0}^{\infty}{\frac{2^{k+1} \cdot 2^{-(k+1)(s+1)}}{k+1}} - \sum_{k=0}^{\infty}{\frac{2^{-(k+1)(s+1)}}{k+1}}\right) = $$
		$$=2^{s+1} \left( -\ln \left( 1-\frac{1}{2^s}\right) - \ln \left( 1-\frac{1}{2^{s+1}}\right)\right)=2^{s+1} \varphi_2(s),$$
		where
		$$\varphi_p(s)=\ln \frac{p^{s+1}-1}{p^{s+1}-p}.$$
		Similarly, in case of $p=3$ we get
		$$F_3(s+1)=\frac{3^{s+1}}{2} \;\varphi_3(s).$$
		
		Since $p^{s+1}-1 \neq p^{s+1}-p$, then $\varphi_p(s) \neq 0$. As $\varphi_p\left( s + \frac{2 \pi i}{\ln p}\right) = \varphi_p(s)$,
		then the maximum and minimum of $|\varphi_p(s)|$ in the strip $\frac{1}{2} \le \Res s \le \frac{11}{10}$ coincide with the maximum and minimum of $|\varphi_p(s)|$
		in rectangular 
		\begin{equation} \label{rectangular}
			\frac{1}{2} \le \Res s \le \frac{11}{10},\quad 0 \le \Ims s \le \frac{2\pi}{\ln p}.
		\end{equation}
		Since $\varphi(s)$ is analytic in (\ref{rectangular}) and $\varphi_p(s) \neq 0$ in (\ref{rectangular}), 
		then, according to the maximum principle (see \cite{rockafellar}, ch. 32) the extremal values $|\varphi_p(s)|$ have to be reached on the boundary of (\ref{rectangular}).
		Since $\varphi_p(s)$ is periodic, then it is enough to examine the behaviour of $s$ along the three sides of (\ref{rectangular}).
		
		By calculations,
		$$\max_{\bfrac{0 \le t \le 2\pi/\ln2}{\sigma = 1/2}}{|\varphi_2(s)|} = \left|\varphi_2\left(\tfrac{1}{2}\right)\right| = 0.79168...; \;\;
					\min_{\bfrac{0 \le t \le 2\pi/\ln2}{\sigma = 1/2}}{|\varphi_2(s)|} = 0.23206...$$
		$$\max_{\bfrac{0 \le t \le 2\pi/\ln2}{\sigma = 11/10}}{|\varphi_2(s)|} = \left|\varphi_2\left(\tfrac{11}{10}\right)\right| = 0.36279...; \;\;
					\min_{\bfrac{0 \le t \le 2\pi/\ln2}{\sigma = 11/10}}{|\varphi_2(s)|} = 0.17323...$$
		$$\max_{\bfrac{1/2 \le \sigma \le 11/10}{t = 0}}{|\varphi_2(s)|} = \left|\varphi_2\left(\tfrac{1}{2}\right)\right| = 0.79168...; \;\;
					\min_{\bfrac{1/2 \le \sigma \le 11/10}{t = 0}}{|\varphi_2(s)|} = 0.36272...$$
		Then, in (\ref{rectangular}) we have
		$$A_1 \le |\varphi_2(s)| \le B_1,\;\;A_1=0.17323...,\;\;B_1=0.79168...$$
		Similarly, for $\varphi_3(s)$ we get
		$$\max_{\bfrac{0 \le t \le 2\pi/\ln3}{\sigma = 1/2}}{|\varphi_3(s)|} = \left|\varphi_3\left(\tfrac{1}{2}\right)\right| = 0.64746...; \;\;
					\min_{\bfrac{0 \le t \le 2\pi/\ln3}{\sigma = 1/2}}{|\varphi_3(s)|} = 0.279736...$$
		$$\max_{\bfrac{0 \le t \le 2\pi/\ln3}{\sigma = 11/10}}{|\varphi_3(s)|} = \left|\varphi_3\left(\tfrac{11}{10}\right)\right| = 0.24985...; \;\;
					\min_{\bfrac{0 \le t \le 2\pi/\ln3}{\sigma = 11/10}}{|\varphi_3(s)|} = 0.16642...$$
		$$\max_{\bfrac{1/2 \le \sigma \le 11/10}{t = 0}}{|\varphi_3(s)|} = \left|\varphi_3\left(\tfrac{1}{2}\right)\right| = 0.64746...; \;\;
					\min_{\bfrac{1/2 \le \sigma \le 11/10}{t = 0}}{|\varphi_3(s)|} = 0.24989...$$
		Consequently, in (\ref{rectangular}) we have
		$$A_2 \le |\varphi_3(s)| \le B_2,\;\;A_2=0.16642...,\;\;B_2=0.64746...$$
		
		Now let us estimate the functions
		$$G_p(s+1)=F_p(s+1) \left( 1-\frac{1}{p^s}\right)^{\frac{1}{2}} \left( 1-\frac{1}{p^{2s}}\right)^{-\frac{1}{24}},\quad \text{for}\quad p=2, 3.$$
		If $p=2$ then
				$$|G_p(s+1)| \le 2^{\sigma+1} \left( 1+\frac{1}{2^\sigma}\right)^{\frac{1}{2}} \left( 1-\frac{1}{2^{2\sigma}}\right)^{-\frac{1}{24}} B_1 = $$
				$$ =\frac{2\sqrt{2^{2\sigma}+2^\sigma}}{\left( 1-\frac{1}{2^{2\sigma}}\right)^{\frac{1}{24}}} B_1 \le
				\frac{2\sqrt{2^{2 \cdot 2}+2^{1.1}}}{\left( 1-\frac{1}{2^{2\cdot 2}}\right)^{\frac{1}{24}}} B_1 < 2.245028 \cdot B_1 < 4.1524;$$
				$$|G_p(s+1)| > 2^{\sigma+1} \left( 1-\frac{1}{2^\sigma}\right)^{\frac{1}{2}} \left( 1+\frac{1}{3^{2\sigma}}\right)^{-\frac{1}{24}} A_1 =$$
				$$=\frac{2 \sqrt{2^{2\sigma}-2^\sigma}}{\left(1+\frac{1}{2^{2\sigma}}\right)^{\frac{1}{24}}} A_1 >
				\frac{2 \sqrt{2-\sqrt{2}}}{\left(\frac{3}{2}\right)^{\frac{1}{24}}} A_1 > 1.505090 \cdot A_1 > 0.260726.$$
		In case $p=3$ we have
				$$|G_p(s+1)| \le \frac{3^{\sigma+1}}{2} \left( 1+\frac{1}{3^\sigma}\right)^{\frac{1}{2}} \left( 1-\frac{1}{3^{2\sigma}}\right)^{-\frac{1}{24}} B_2 =$$
				$$=\frac{\frac{3}{2}\sqrt{3^{2\sigma}+3^\sigma}}{\left( 1-\frac{1}{3^{2\sigma}}\right)^{\frac{1}{24}}} B_2 \le
				\frac{\frac{3}{2}\sqrt{3^{2 \cdot 2}+3^{1.1}}}{\left( 1-\frac{1}{3^{2\cdot 2}}\right)^{\frac{1}{24}}} B_2 < 5.745949  B_2 < 3.720278;$$
				$$|G_p(s+1)| > \frac{3^{\sigma+1}}{2} \left( 1-\frac{1}{3^\sigma}\right)^{\frac{1}{2}} \left( 1+\frac{1}{3^{2\sigma}}\right)^{-\frac{1}{24}} A_2 =$$
				$$=\frac{\frac{3}{2} \sqrt{3^{2\sigma}-3^\sigma}}{\left(1+\frac{1}{3^{2\sigma}}\right)^{\frac{1}{24}}} A_2 >
				\frac{\frac{3}{2} \sqrt{3-\sqrt{3}}}{\left(\frac{4}{3}\right)^{\frac{1}{24}}} \cdot A_2 > 1.668923 \cdot A_2 > 0.277752.$$
	Thus,
	$$|G(s+1)| = \prod_p{|G_p(s+1)|} \le |G_2(s+1)| \cdot |G_3(s+1)| \cdot \exp \left({\sum_p{\ln |G_p(s+1)|}} \right) \le $$
	$$\le	4.1524 \cdot 3.720278 \cdot \exp \left({\sum_{p \ge 5}{\frac{C_7}{p^{3/2}}}} \right) <
	15.4481 \cdot e^{0.1605 \cdot C_7} < 15.4481 \cdot e^{0.52} < 26,$$
	and
	$$|G(s+1)| > 0.260726 \cdot 0.277752  \cdot e^{-0.52} > 0.04305 > \frac{1}{24}.$$
	
	Finally we have $\frac{1}{24} < |G(s)| < 26$ in the strip $\frac{3}{2} \le \sigma \le \frac{21}{10}$,.
	
	Now let us consider the integral $I$.
	Let $\Gamma$ be the boundary of the rectangle with the vertices $\frac{3}{2} \pm iT, b \pm iT$, $ 2 \le T \le x$, 
	where the zeros of $\zeta(s-1)$ of the form $\rho = \frac{3}{2} + i\gamma$, $|\gamma| < T$, are avoided 
	by the semicircles of the infinitely small radius lying to the right of the line $\Res s = \frac{3}{2}$,
	the pole of $\zeta(s-1)$ at the point $s=2$ is avoided be the arcs $\Gamma_1, \Gamma_2$ with the radius $\frac{1}{\ln x}$, 
	and let a horizontal cut be drawn from the critical line inside this rectangle to each zero $\rho = \beta +i\gamma, \frac{3}{2} < \beta < 2, |\gamma| <T$. 
	Then the function $|F(s)|$ is analytic inside $\Gamma$.
		By the Cauchy residue theorem,
		$$j_0 = - \sum_{k=1}^{8}{j_k} - \sum_{\rho}{j_{\rho}} = - (j_4+j_5) - \sum_{k \neq 4, 5}{j_k} - \sum_{\rho}{j_{\rho}}.$$
	
	\begin{figure}[tbh]
			\begin{center}
				\includegraphics{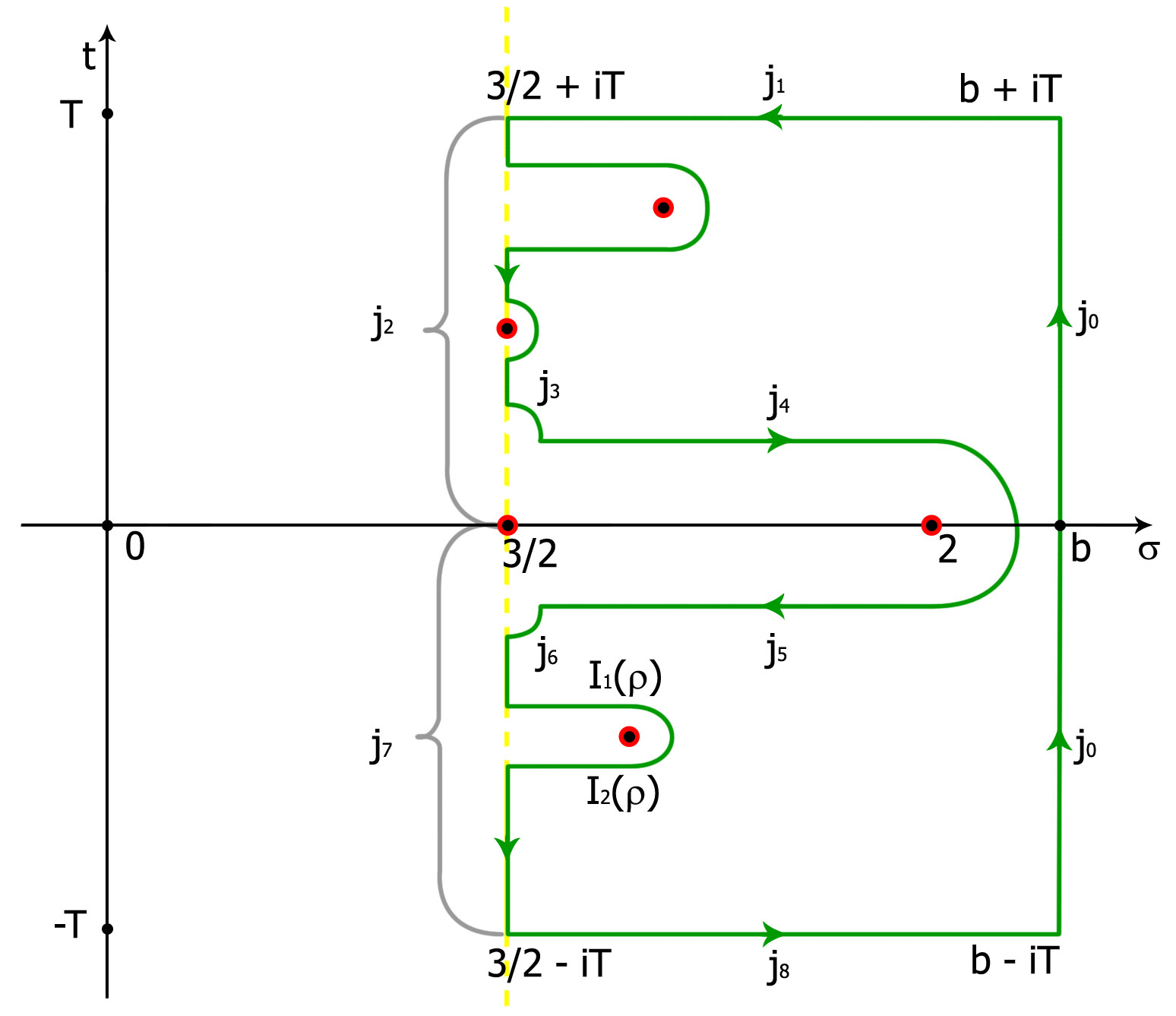}
			\end{center}
	\end{figure}
	By the proof of theorem 1, we have on $\Gamma$:
	$$F(s) = \frac{(\zeta(s-1))^{\frac{1}{2}}}{{\zeta(2(s-1))}^{\frac{1}{24}}} \cdot G(s).$$
	According to the lemma \ref{T^()lnT},
	$$|F(s)| \ll T^{\frac{1-(\sigma-1)}{6}}(\ln T)^{\frac{13}{24}} = T^{\frac{2-\sigma}{6}} (\ln T)^{\frac{13}{24}}.$$
	Then
	$$|j_1| \ll \frac{1}{T} \int_{\frac{3}{2}}^{-b}{x^{\sigma} T^{\frac{2-\sigma}{6}} (\ln T)^{\frac{13}{24}} d\sigma} \ll \frac{x^2}{T} \int_{\frac{3}{2}}^{-b}{\left(\frac{x}{T^{\frac{1}{6}}}\right)^{\sigma-2}(\ln T)^{\frac{13}{24}} d\sigma} \ll \frac{x^2}{T} (\ln T)^{\frac{13}{24}}.$$
	The similar estimate is valid for $j_8$.
	
	Next, on the arcs $\Gamma_1$ and $\Gamma_2$ we have
	
	\begin{figure}[tbh]
			\begin{center}
				\includegraphics{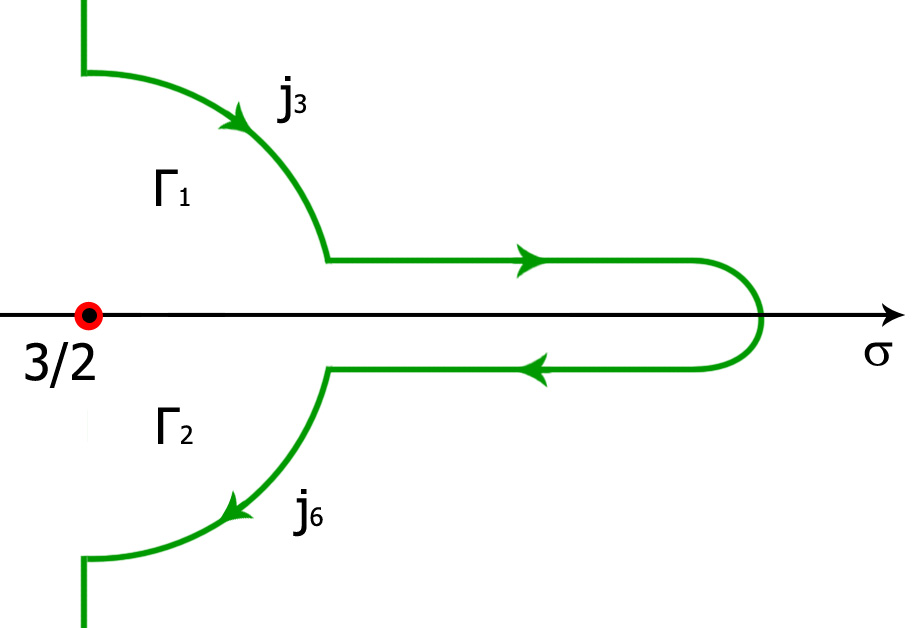}
			\end{center}
	\end{figure}
	
	$$|\zeta(s-1)| = \left| \frac{1}{2} + \frac{1}{s-2} + (s-1)\int_{1}^{\infty}{\frac{\rho(u)}{u^s}} du\right| \le 0.5 +1+0.1 + 1.6 \cdot 0.5 \int_{1}^{\infty}{\frac{du}{u^{\frac{3}{2}}}} < 3.2,$$
	$$|\zeta(2(s-1))| \ge \frac{1}{|2s-3|} -1.1 \ge 0.5\ln x - 1.1 > 0.4 \ln x.$$
	
	Hence
	$$|F(s)| \le \frac {(3.2)^{\frac{1}{2}}}{{(0.4 \ln x)}^{\frac{1}{24}}} < 1.$$
	Therefore,
	$$|j_3+j_6| \le \frac{1}{2\pi} \int_{-\frac{\pi}{2}}^{\frac{\pi}{2}}{\frac{2(2x)^{\frac{3}{2} + \frac{1}{\ln x}}}{\frac{5}{2}} \frac{d \varphi}{\ln x}} \ll \frac{x^{\frac{3}{2}}}{\ln x}.$$
	Since
	$$|F(s)| \le |\zeta(s-1)|^{\frac{1}{2}}(\ln x)^{\frac{1}{24}}|G(s)| \ll (\ln x)^{\frac{1}{24}} \left| (\zeta(\sigma-1+it)) \right| ^{\frac{1}{2}},$$
	then
	$$|j_2| = \left| \text{p.v.} \frac{1}{2 \pi i}\int_{\frac{3}{2} + iT}^{\frac{3}{2} + \frac{i}{\ln x}}{F(s) \cdot \frac {(x + h)^s - x^s}{s} ds} \right|
			\ll \int_{\frac{1}{\ln x}}^{T} {(\ln x)^{\frac{1}{24}} \cdot |\zeta(\sigma -1 + it)|^{\frac{1}{2}}} x^{\frac{3}{2}} \frac{dt}{t+1} \ll $$
			$$\ll (\ln x)^{\frac{1}{24}} x^{\frac{3}{2}} \int_{0}^{T}{\frac{\left|\zeta\left(\frac{1}{2}+it\right)\right|^{\frac{1}{2}} dt}{t+1}} =
			(\ln x)^{\frac{1}{24}} x^{\frac{3}{2}} \sum_{\nu \ge 0}{\int_{T/{2^\nu}}^{T/2^{\nu +1}}{\frac{\left|\zeta\left(\frac{1}{2} +it\right)\right|^{\frac{1}{2}}}{t+1}dt}}.$$
	
	Denoting the summands in the last sum by $j(\nu)$ and taking $X=T\cdot2^{-\nu}$, by Hölder inequality we get:
		$$j(\nu) \ll \frac{1}{X} \left( \int_{X}^{2X}{\left|\zeta\left(\tfrac{1}{2}+it\right)\right|^2 dt}\right)^{\frac{1}{4}} X^{1-\frac{1}{4}} \ll
		\frac{1}{X} \left( X\ln X\right)^{\frac{1}{4}}X^{\frac{3}{4}} \ll
		(\ln X)^{\frac{1}{4}} \ll (\ln T)^{\frac{1}{4}}.$$
	Thus,
	$$\sum_{\nu \ge 0}{j(\nu)} \ll (\ln T)^{\frac{5}{4}}.$$
	Then
	$$|j_2| \ll (\ln x)^{\frac{1}{24}} x^{\frac{3}{2}} (\ln x)^{\frac{5}{4}} \ll (\ln x)^{\frac{31}{24}} x^{\frac{3}{2}}.$$
	The integral $j_7$ is estimated as above.
		
	The main term arises from $j_4$ and $j_5$.
	Let us define the entire function $w(s)$ by the relation
	$$\zeta(s-1) = \frac{w(s-1)}{s-2}.$$
	Let $s = 2-u + i \cdot 0$, where $0 \le u \le \frac{1}{2}$. Then
	$$\sqrt{\zeta(s-1)} = \frac {\sqrt{w(s-1)}}{\sqrt{-u+i \cdot 0}}.$$
	Since $-u+i\varepsilon \to u\cdot e^{\pi i}$ as $\varepsilon \to +0$ then
	$$\sqrt{-u+i \cdot 0} = \sqrt{u} e^{\frac{\pi i}{2}}, \;\; \sqrt{\zeta(s-1)}=\frac{\sqrt{w(\sigma-1)}}{\sqrt{u}} e^{-\frac{\pi i}{2}}.$$
	Therefore, on the upper edge of the cut we have
		$$F(s) = \frac{\sqrt{w(1-u)}}{(\zeta(2-2u))^{\frac{1}{24}}} \;G(2-u)\;\frac{e^{-\frac{\pi i}{2}}}{\sqrt{u}} = \frac{\Pi(u)e^{-\frac{\pi i}{2}}}{\sqrt{u}} = \frac{\Pi(u)}{i \sqrt{u}},$$
		where
		$$\Pi(u) = G(2-u) \frac{\sqrt{w(1-u)}}{(\zeta(2-2u))^{\frac{1}{24}}}.$$
		
	Then for $j_4$ we obtain:
		$$j_4 = \frac{1}{2 \pi i} \int_{\frac{3}{2} + \frac{1}{\ln x} + i \cdot 0}^{2+i \cdot 0}{F(\sigma + i \cdot 0) \frac{(x+h)^s- x^s}{s} ds} = $$
		$$=\frac{1}{2\pi i} \int_{\frac{1}{2}}^{0}{\frac{\Pi(u)}{i \sqrt{u}} \frac{(x+h)^{2-u}-x^{2-u}}{2-u}} (-du) = $$
		$$= -\frac{1}{2\pi} \int_{0}^{\frac{1}{2}}{\frac{\Pi(u)}{\sqrt{u}} \int_{x}^{x+h}{y^{1-u}dy}du} =$$
		$$= \frac{1}{2\pi} \int_{x}^{x+h}{y \int_{0}^{\frac{1}{2}}{\frac{\Pi(u)}{\sqrt{u}} y^{-u} du}dy}.$$
		Suppose that $N \ge 0$ is fixed. Then
		$$\Pi(u)= \Pi_0+\Pi_1u + \Pi_2u^2 + \ldots + \Pi_N u^N + O_N(u^{N+1})$$
		and
		$$j_4 = -\frac{1}{2\pi} \int_{x}^{x+h} {y \left(\sum_{0 \le n \le N}{\Pi_n \int_{0}^{\frac{1}{2}}{\frac{u^{n}y^{-u}}{\sqrt{u}}du}}\right)dy} + O\left( J\right)=$$
		$$=-\frac{1}{2\pi} \sum_{0 \le n \le N}{\Pi_n \Gamma\left( n + \frac{1}{2}\right) \int_{x}^{x+h}{\frac{y \; dy}{(\ln y)^{n+\frac{1}{2}}}}} + O\left( J\right),$$
		where
		$$J =  -\frac{1}{2\pi}\int_{x}^{x+h} {y {\int_{0}^{\frac{1}{2}}{\frac{u^{N+1}y^{-u}}{\sqrt{u}}}du}\;dy} \ll \frac{xh}{(\ln x)^{N+\frac{3}{2}}}{\Gamma\left(N+\frac{3}{2}\right)}.$$
		
		Let us evaluate a contribution of $n$-th term to the sum. Assume
		$$\varphi(y) = \frac{y}{(\ln y)^{\nu}}, \;\; \varphi{'}(y) = \frac{y \ln y - \nu}{y (\ln y)^{\nu+1}}.$$
		If $x \le y \le x+h$ then
		$$\varphi(y) = \varphi(x) + (y-x) \; \theta \; \frac{x \ln x + \nu}{x(\ln x)^{\nu+1}},\;\;|\theta| \le 1.$$
		Hence,
		$$\int_{x}^{x+h}{\varphi(y) dy} = h \varphi(x) + \theta_1 \frac{x \ln x + \nu}{x(\ln x)^{\nu +1}} \int_{x}^{x+h}{(y-x)dy} = 
		\frac{xh}{(\ln x)^\nu} + \frac{\theta_1}{2} \;h^2\; \frac {x \ln x + \nu}{x(\ln x)^{\nu +1}}.$$
		Finally, we state: 
		$$j_4 = -\frac{xh}{2\pi} \sum_{0 \le n \le N}{\frac{\Pi_n \Gamma\left( n+\frac{1}{2}\right)}{(\ln x)^{n+\frac{1}{2}}}}+ O\left( \frac{h}{(\ln x)^{N+\frac{3}{2}}}\right) + O\left( \frac{h^2}{x} \frac{1}{(\ln x)^{N+\frac{3}{2}}}\right)$$
		and
		$$-(j_4+j_5)=
		-\frac{xh}{2 \pi (\ln x)^{\frac{1}{2}}} \left( \sum_{0 \le n \le N}{\frac{\Pi_n \Gamma \left( n+\frac{1}{2}\right)}{(\ln x)^n}} + O\left( \frac{1}{(\ln x)^{N+1}}\right) + O \left( \frac{h}{x(\ln x)^{N+1}}\right)\right).$$
		
		Now let us estimate the sum
		$$\sum_{|\gamma|<T}{j_{\rho}},\quad{\rm where}\quad j_{\rho}=I_1(\rho)+I_2(\rho).$$
		
		\begin{figure}[tbh]
			\begin{center}
				\includegraphics{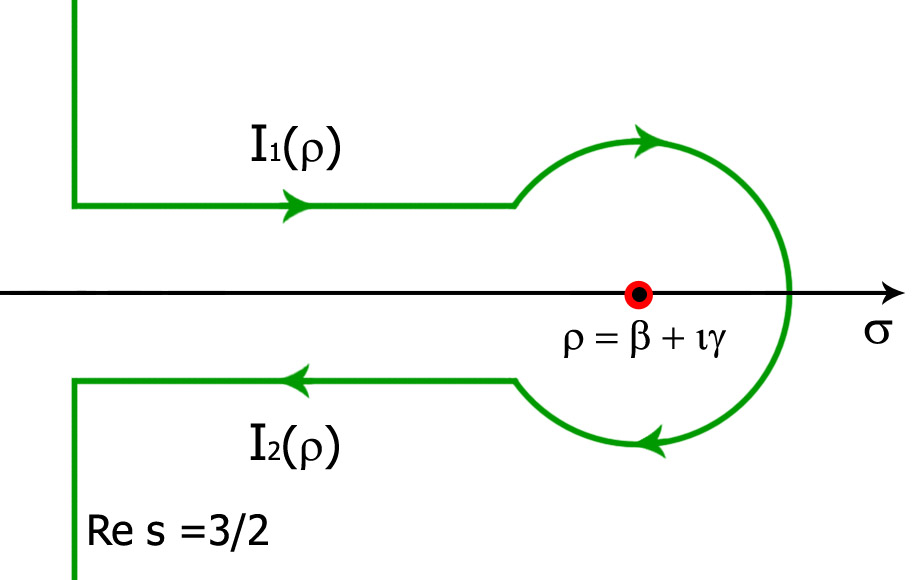}
			\end{center}
		\end{figure}
		By analogy with the proof of the theorem \ref{theorem1}, we get
		$$I_1(\rho) =  \int_{\frac{3}{2}+ i \gamma}^{\beta+1+ i \gamma}{F(s) \frac{(x+h)^s - x^s}{s} ds} \ll \int_{\frac{3}{2}}^{\beta+1}{\left|F(s)\right| h x^{\sigma-1} d\sigma} \ll$$
		$$\ll h \int_{\frac{3}{2}}^{\beta +1}{T^{\frac{2-\sigma}{6}} (\ln x)^{\frac{13}{24}} x^{\sigma-1} d\sigma}
		\ll xh (\ln x)^{\frac{13}{24}} \int_{\frac{3}{2}}^{\beta+1}{ \left( \frac{T^{\frac{1}{6}}}{x}\right)^{2-\sigma} d\sigma}.$$
		Hence,
		$$|j_\rho| \ll	x h (\ln x)^{\frac{13}{24}} \int_{\frac{3}{2}}^{\beta+1}{\left( \frac{T^{\frac{1}{6}}}{x}\right)^{2-\sigma} d\sigma}
		\ll x h (\ln x)^{\frac{13}{24}} \int_{\frac{1}{2}}^{\beta}{\left( \frac{T^{\frac{1}{6}}}{x}\right)^{1-v} dv} = $$
		$$ = x h (\ln x)^{\frac{13}{24}} \int_{\frac{1}{2}}^{1}{ g(\rho; v) \left( \frac{T^{\frac{1}{6}}}{x}\right)^{1-v} dv},$$
		where
		$$g(\rho; v) = \left\{
			\begin{aligned}
				1,\;\text{if}\;v \le \beta,\\ 
				0,\;\text{if}\;v > \beta.\\
			\end{aligned}
			\right.$$
		Applying the lemma \ref{density}, we obtain:
		$$\sum_{|\gamma|<T}{j_\rho} \ll x h (\ln x)^{\frac{13}{24}} \int_{\frac{1}{2}}^{1-\varrho (T)} { \left( \sum_{|\gamma| < T}{g(\rho; v)}\right)\left( \frac{T^{\frac{1}{6}}}{x}\right)^{1-v} dv} \ll$$
		$$\ll  x h (\ln x)^{\frac{13}{24}} \int_{\frac{1}{2}}^{1-\varrho (T)} {N(v; T) \left( \frac{T^{\frac{1}{6}}}{x}\right)^{1-v} dv} \ll
		x h (\ln x)^{44+\frac{13}{24}} \left( \frac{T^{\frac{12}{5}+\frac{1}{6}}}{x}\right)^{1-\varrho (T)} \ll$$
		$$\ll x h (\ln x)^{45} \left( \frac{T^{\frac{77}{30}}}{x}\right)^{1-\varrho (T)}.$$
		Choosing $T$ from the equation
		$$T^{\frac{77}{30}} = x \cdot D^{-1}(x), \;D(x)=e^{C_2 (\ln x)^{0.8}},$$
		we get
		$$T= x^{\frac{30}{77}} D(x)^{-\frac{30}{77}}.$$
		Now we can easily conclude that the formula
		$$j_0 = \frac{xh}{2 \pi(\ln x)^{\frac{1}{2}}} \left( \sum_{0 \le n \le N}{\frac{\Pi_n \Gamma(n+\frac{1}{2})}{(\ln x)^n}}\right) + O(J),$$
		where
		$$J = \frac{x^2}{T} (\ln T)^{\frac{13}{24}} + x^{\frac{3}{2}}(\ln x)^{\frac{31}{24}} + \frac{x^{\frac{3}{2}}}{\ln x} + \frac{xh}{(\ln x)^{N+\frac{3}{2}}}+ \frac{h^2}{(\ln x)^{N+\frac{3}{2}}} + xh(\ln x)^{45} \left( \frac{T^{\frac{77}{30}}}{x}\right)^{1-\varrho (T)} \ll$$
		$$\ll \frac{x^2}{T} (\ln T)^{\frac{5}{6}} + \frac{xh}{(\ln x)^{N+\frac{3}{2}}}+ \frac{h^2}{(\ln x)^{N+\frac{3}{2}}} + xh(\ln x)^{45} \left( \frac{T^{\frac{77}{30}}}{x}\right)^{1-\varrho (T)},$$
		is an asymptotic, if
		$$h = x^{\alpha} e^{C_2 (\ln x)^{0.8}} \gg \frac{x}{T} (\ln x)^2,$$
		where
		$$\alpha=1 - \frac{30}{77} = \frac{47}{77}.$$
	\end{subsection}

	\begin{subsection}{The mean-value of the function $\frac{\D 1}{\D r(n)}$ on the short interval}
		Suppose that $\sigma = \Res s >1$ and let
		$$\varkappa(n) = \displaystyle{\sum_{d|n}{\chi_4(d)}} = \displaystyle{\frac{1}{4}}\;r(n),\;\;\;F(s) = \displaystyle{\sum_{n=1}^{\infty}{'}{\frac{1}{\varkappa(n)} \cdot n^{-s}}},$$
		where $\sum{'}$ means that the summing is going over all $n$, with $\varkappa(n) \neq 0$.
		This series converges absolutely, since
		$$|F(s)| \le \sum_{n=1}^{\infty}{'}{|a_n| \cdot n^{-\sigma}} \le \sum_{n=1}^{\infty}{'}{n^{-\sigma}} \le 1 + \int_{1}^{\infty}{\frac{du}{u^\sigma}} = 1 + \frac{1}{\sigma-1}.$$
		Setting $a_n=\displaystyle {\frac{1}{\varkappa(n)}}$, $A(n) \equiv 1$, $b = 1 + \displaystyle{\frac{1}{\ln x}}$, $\alpha=1$ in the lemma \ref{perron}, we get:
		$$S_3 = S(x,h;f_3) = I + O(R),$$
		where
		$$I=\frac{1}{2\pi i} \int_{b-iT}^{b+iT}{F(s) \frac{x^s}{s} ds},\;R=\frac{x^b}{T(b-1)} + \frac{x A(2x) \ln x}{T} \ll \frac{x \ln x}{T}.$$
		Further, let $F(s)=\prod_{p}{'}F_p(s)$.
		Then
			$$F_2(s) = 1 + \frac{1}{2^s}+\frac{1}{2^{2s}} + \ldots = \frac{1}{1-\frac{1}{2^s}} = \left( 1 - \frac{1}{2^s}\right)^{-1}.$$
		In the case $p \equiv 1 (\mmod 4)$, we have
			$$F_p(s) = 1 + \frac{1}{2p^s}+\frac{1}{3p^{2s}} + \ldots = \left(1-\frac{1}{p^s}\right)^{-\frac{1}{2}}\left(1-\frac{1}{p^{2s}}\right)^{\frac{1}{24}} G_p(s),$$
			where the function $G_p(s)$ is defined in the proof ot the theorem \ref{theorem1}.
		Finally, for $p \equiv 3 (\mmod 4)$ it is true that
			$$F_p(s) = 1 + \frac{1}{p^{2s}}+\frac{1}{p^{4s}}+\ldots = \left(1-\frac{1}{p^{2s}}\right)^{-1}.$$
		
		Let us write $F(s)$ in the form
		$$F(s)=\left(1-\frac{1}{2^s}\right)^{-1} \prod_{p \equiv 1(\mmod 4)}{\left(1-\frac{1}{p^s}\right)^{-\frac{1}{2}}\left(1-\frac{1}{p^{2s}}\right)^{\frac{1}{24}} G(s)} \prod_{p \equiv 3(\mmod 4)}{\left(1-\frac{1}{p^{2s}}\right)^{-1}},$$
		where
		$$G(s)=\prod_p{G_p(s)}=\prod_{p \equiv 1(\mmod 4)}{\left(1-\frac{1}{p^s}\right)^{\frac{1}{2}}\left(1-\frac{1}{p^{2s}}\right)^{-\frac{1}{24}}}\left( 1+\frac{1}{2p^s} +u(s)\right),$$
		$$u(s)=\frac{1}{3p^{2s}}+ \ldots$$
		
		Applying the arguments used in the proof of Theorem 1 in the case $k=2$, we get
		$$|\ln G_p(s)| \le \frac{19}{2} \cdot \frac{1}{p^{3\sigma}} \le \frac{19}{2p^{\frac{3}{2}}}.$$
		
		Thus,
		$$|\ln G(s)| \le \frac{19}{2} \sum_{p \equiv 1(\mmod 4)}{\frac{1}{p^{\frac{3}{2}}}} \le \frac{19}{2} \sum_{p \ge 5}{\frac{1}{p^{\frac{3}{2}}}} < 1.525,$$
		
		$$e^{-1.525} < |G(s)| < e^{1.525},\;\;0.2 <|G(s)|<4.6.$$
		
		Let
		$$f(s)=\prod_{p \equiv 1(\mmod 4)}{\left( 1-\frac{1}{p^s}\right)^{-1}}.$$
		Since
		$$\frac{1+\chi_4(p)}{2} = 
		\left\{
			\begin{aligned}
				1,\;\text{if}\;p \equiv 1(\mmod 4),\\ 
				0,\;\text{if}\;p \equiv 3(\mmod 4),\\
				\frac{1}{2},\;{\rm if}\;p\;{\rm is}\;{\rm an}\;{\rm odd}\;{\rm number,}\\
			\end{aligned}
			\right.$$
		then
		$$\ln f(s)=-\sum_{p \equiv 1(\mmod 4)}{\ln \left( 1-\frac{1}{p^s}\right)} = -\sum_{p}{\frac{1+\chi_4(p)}{2} \ln\left( 1-\frac{1}{p^s}\right)} + \frac{1}{2}\ln \left( 1-\frac{1}{2^s}\right)=$$
		$$=\frac{1}{2} \ln \zeta(s) + \frac{v(s)}{2} + \frac{1}{2} \ln \left( 1-\frac{1}{2^s}\right),$$
		where
		$$v(s)=-\sum_{p}{\chi_4(p) \ln \left( 1-\frac{1}{p^s}\right)}.$$
		Using the relation
		$$\ln L(s,\chi_4) = -\sum_{p}{\ln \left( 1-\frac{\chi_4(p)}{p^s}\right)},$$
		we obtain
		$$v(s) = \sum_{p}{\left( \ln \left( 1- \frac{\chi_4(p)}{p^s}\right) - \chi_4(p)\ln \left( 1-\frac{1}{p^s}\right) - \ln \left(1-\frac{\chi_4(p)}{p^s}\right)\right)} =
		\ln L(s,\chi_4) + w(s),$$
		where
		$$w(s) = \sum_{p}{w_p(s)},\;\;w_{p}(s)=\ln \left( 1-\frac{\chi_4(p)}{p^s}\right) - \chi_4(p)\ln \left( 1-\frac{1}{p^s}\right),$$
		with
		$$w_p(s) = 
		\left\{
			\begin{aligned}
				0,\;{\rm if}\;p=2, p \equiv 1(\mmod 4),\\ 
				\ln \left( 1-\displaystyle{\frac{1}{p^{2s}}}\right),\;{\rm if}\;p \equiv 3(\mmod 4).\\
			\end{aligned}
			\right.$$
		Then for $v(s)$ it holds that:
		$$v(s) = \ln L(s, \chi_4) + \sum_{p \equiv 3(\mmod 4)}{\ln \left( 1-\frac{1}{p^{2s}}\right)}.$$
		Hence,
		$$\ln f(s) = \frac{1}{2} \ln \zeta(s) + \frac{1}{2} \ln L(s, \chi_4) + \frac{1}{2}\sum_{p \equiv 3(\mmod 4)}{\ln \left( 1-\frac{1}{p^{2s}}\right)} + \frac{1}{2} \ln \left( 1-\frac{1}{2^s}\right),$$
		Thus,
		$$F(s)= \left( 1-\frac{1}{2^s}\right)^{-1} \left( f(s)\right)^{\frac{1}{2}} \left( f(2s)\right)^{-\frac{1}{24}} G(s) \prod_{p \equiv 3(\mmod 4)}{\left( 1-\frac{1}{p^{2s}}\right)^{-1}}=$$
		$$= \left( 1-\frac{1}{2^s}\right)^{-1} \left( \zeta(s) L(s, \chi_4)\right)^{\frac{1}{4}} \prod_{p \equiv 3(\mmod 4)} {\left( 1-\frac{1}{p^{2s}}\right)^{-\frac{3}{4}}} \left( 1-\frac{1}{2^s}\right)^{\frac{1}{4}} G(s) \left(f(2s) \right)^{-\frac{1}{24}}=$$
		$$= \left( 1-\frac{1}{2^s}\right)^{-\frac{3}{4}} \left( \zeta(s) L(s, \chi_4)\right)^{\frac{1}{4}} \left( \zeta(2s) \right)^{\frac{3}{4}} \left( 1-\frac{1}{2^{2s}}\right)^{\frac{3}{4}} \left( f(2s)\right)^{-\frac{19}{24}} G(s)=$$
		$$= (\zeta(s) L(s,\chi_4))^{\frac{1}{4}}  (\zeta(2s))^{\frac{17}{48}}  (L(2s,\chi_4))^{-\frac{19}{48}} H(s),$$
		where
		$$H(s)=\left(1-\frac{1}{2^s}\right)^{-\frac{3}{4}}  \left(1 - \frac{1}{2^{2s}}\right)^{\frac{17}{48}}  \prod_{p \equiv 3(\mmod 4)}{\left(1 - \frac{1}{p^{4s}}\right)^{-\frac{19}{48}}}  G(s).$$
		For $H(s)$ we have:
		$$\ln H(s) = -\frac{3}{4} \ln \left( 1-\frac{1}{2^s}\right) + \frac{17}{48} \ln \left( 1 - \frac{1}{2^{2s}}\right) + \ln G(s) -\frac{19}{48} \sum_{p \equiv 3(\mmod 4)}{\ln H_p(s)},$$
		where
		$$H_p(s)=1-\frac{1}{p^{4s}}.$$
		Then
		$$|\ln H_p(s)| \le \sum_{k=0}^{\infty}{\left| \frac{1}{p^{4s(k+1)} (k+1)}\right|} \le \frac{1}{p^{4\sigma}} \frac{1}{1-\frac{1}{p^{4\sigma}}} \le \frac{1}{p^{4\sigma}} \frac{1}{1-\frac{1}{p^2}} \le \frac{2}{p^{4\sigma}}.$$
		Next
		$$|\ln H(s)| \le C_1 + \frac{19}{24} \sum_{p \equiv 3(\mmod 4)}{\frac{1}{p^{4\sigma}}} \le C_2,$$
		so that
		$$e^{-C_2} \le |H(s)| \le e^{C_2}.$$
		
		Suppose that $2 \le T \le x$ does not coincide with the ordinate of zeros of $\zeta(s)$ and $L(s, \chi_4)$.
		Let $\Gamma$ be the boundary of the rectangle with the vertices $\frac{1}{2} \pm iT, b \pm iT$,
		where the zeros of $\zeta(s)$ and $L(s,\chi_4)$ are avoided by semicircles of the infinitely small radius lying to the right of the line $\Res s = \frac{1}{2}$,
		the point $s = \frac{1}{2}$ is avoided by two arcs $\Gamma_1, \Gamma_2$ with the radius $\frac{1}{\ln x}$,
		and let a horizontal cut be drawn from the critical line inside this rectangle to each zero $\rho = \frac{1}{2} + i\gamma,  \frac{1}{2} < \beta < 1,|\gamma| <T$.
		Then the function $F(s)$ is analytic inside $\Gamma$.
		Applying the Cauchy residue theorem,
		
		$$j_0 = - \sum_{k=1}^{8}{j_k} - \sum_{\rho}{j_{\rho}} = - (j_4+j_5) - \sum_{k \neq 4, 5}{j_k} - \sum_{\rho}{j_{\rho}}.$$
		
		\begin{figure}[tbh]
			\begin{center}
				\includegraphics{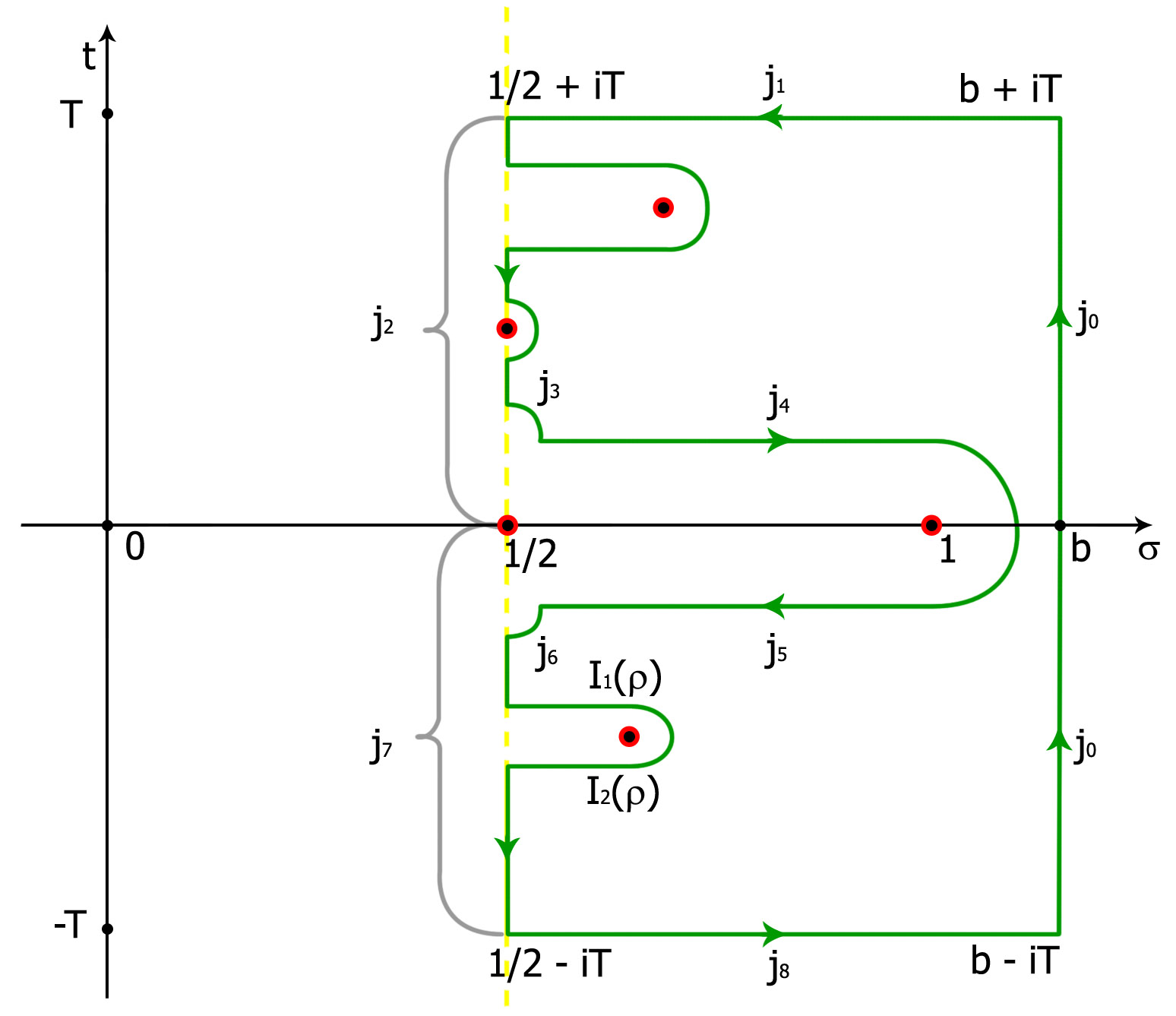}
			\end{center}
		\end{figure}
		
		Consider the integrals $j_1$ and $j_8$.
		For $\Res s \ge \frac{1}{2}+\frac{1}{\ln x}$ we have:
		$$|\zeta(2s)| \le 1 + \sum_{n=2}^{\infty}{\frac{1}{n^{2\sigma}}} \le 1 + \frac{1}{2\sigma-1} \le 1+ \frac{\ln x}{2} < \ln x.$$
		
		Suppose that $|s-1| \le \frac{1}{\ln x}$. If $\Res s >1$, then
		$$L(s,\chi_4) = \sum_{n=1}^{\infty}{\frac{\chi_4(n)}{n^s}}.$$
		This series converges absolutely. Then,
		$$L'(s,\chi_4) = -\sum_{n=2}^{\infty}{\frac{\chi_4(n) \ln n}{n^s}}.$$
		By Abel summation formula,
		\begin{equation} \label{abel}
			L^{'}(s,\chi_4) = -s \int_{1}^{\infty}{\frac{\mathbb{C}(u)}{u^{s+1}} du,}
		\end{equation}
		where 
		$$\mathbb{C}(u)=\sum_{1<n\le u}{\chi_4(u) \ln u}.$$
		The integral in (\ref{abel}) converges for $\Res s >0$. It means that we can use (\ref{abel}) for $|s-1| \le \frac{1}{\ln x}$.
		Let us note that
		$$\left| \mathbb{C}(u) \right| < \ln u,$$
		uniformly in $u \ge 1$.
		Indeed, if $u=4m+1$, then
		$$\mathbb{C}(u)=\ln (4m+1)-(\ln(4m-1)-\ln(4m-3))-(\ln(4m-5)-\ln(4m-7))-\ldots < \ln (4m+1) < \ln u,$$
		$$\mathbb{C}(u)=(\ln(4m+1)-\ln(4m-1))+(\ln(4m-3)-\ln(4m-5))+\ldots > 0,$$
		consequently
		$$0< \mathbb{C}(u) < \ln u$$
		(the cases $u \equiv 0, 2, 3 (\mmod 4)$ are treated as above).
		Since $\mathbb{C}(u)=\mathbb{C}([u])$, then for $|s-1| \le \delta$ we have
		$$|L(s, \chi_4)| \le (1+\delta) \int_{1}^{\infty}{\frac{\ln u}{u^{2-\sigma}}du} \le \frac{1+\delta}{(1-\delta)^2}.$$
		Next,
		$$L(s, \chi_4) = L(1, \chi_4) + \int_{1}^{s}{L^{'}(u, \chi_4) du},$$
		and
		$$|L(s, \chi_4)| \ge |L(1, \chi_4)| - |s-1|\frac{1+\delta}{(1-\delta)^2}.$$
		Since
		$$L(1, \chi_4)=1 - \frac{1}{3} + \frac{1}{5} - \frac{1}{7}+ \ldots = \frac{\pi}{4},$$
		we have
		$$|L(s, \chi_4)| \ge \frac{\pi}{4}- |s-1|\frac{1+\delta}{(1-\delta)^2}.$$
		Taking $\delta = \frac{1}{4}$, for $|s-1| \le \frac{1}{\ln x}$, we obtain
		
		\begin{equation} \label{evofLag}
			|L(s, \chi_4)| \ge \frac{\pi}{4} - \frac{5}{9} > \frac{1}{5}.
		\end{equation}
		
		Thus, if $\frac{1}{2} \le \Res s \le \frac{5}{8}$, then $|L(2s, \chi_4)|^{-1} \le 5$.
		In case of $\frac{5}{8} \le \Res s \le 1$ we have
		$$|L(2s, \chi_4)|^{-1} \le 1 + \sum_{n=2}^{\infty}{\frac{1}{n^{2\sigma}}} \le 1 + \frac{1}{2\sigma - 1} \le 1 + \frac{1}{2 \cdot \frac{5}{8} - 1} \le 5.$$
		Then
		$$F(s) \ll \left( T^{\frac{1-\sigma}{3}} \ln T \right)^2,$$
		Consequently,
		$$|j_1| \ll \frac{1}{T} \int_{\frac{1}{2}}^{b}{T^{\frac{2(1-\sigma)}{3}} (\ln x)^{2} x^\sigma d\sigma} = 
		\frac{x}{T} (\ln x)^{2}\int_{\frac{1}{2}}^{b}{\left( \frac{T^{\frac{2}{3}}}{x}\right)^{1-\sigma} d\sigma} \ll \frac{x}{T} (\ln x)^{2}.$$
		The same estimation is valid for $j_8$.
		
		\begin{figure}[tbh]
			\begin{center}
				\includegraphics{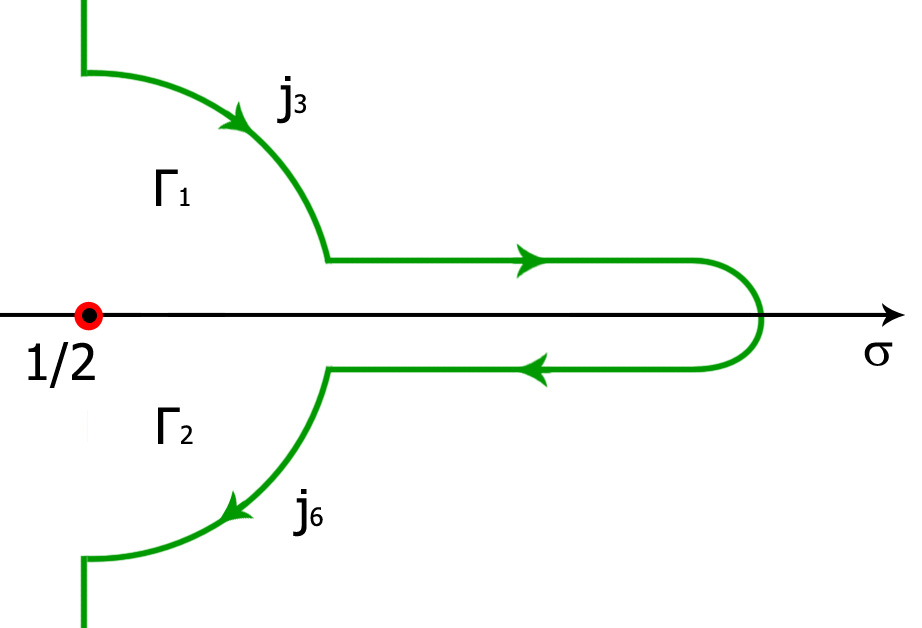}
			\end{center}
		\end{figure}
		By lemma \ref{evnofzeta}, on the arcs $\Gamma_1$, $\Gamma_2$ we have
		$$|\zeta(s)| \le 3.2,\;\;|\zeta(2s)| \le 2.2,\;\;|L(s,\chi_4)| \le 3.2, |\zeta(s)| \le 3.2.$$
		
		Using (\ref{evofLag}), we get
		$$|F(s)| \le (3.2)^{\frac{1}{2}} (2.2)^{\frac{17}{48}} 5^{\frac{19}{48}}  \left| H(s) \right| < C.$$
		Repeating the proofs of the theorems \ref{theorem1} and \ref{theorem2}, we get
		$$|j_3+j_6| \le \frac{C}{2\pi} \int_{\Gamma_1 \cup \Gamma_2}{\left| \frac{(x+h)^s-x^s}{s}\right| ds} \ll \frac{\sqrt{x}}{\ln x}.$$
		
		Since
		$$F(s) \ll \left | (\zeta(s)) \right|^{\frac{1}{2}} (\ln T)^{\frac{1}{2}},\quad{\rm then}$$
		$$|j_2| \ll \left| \int_{\frac{1}{\ln x}}^{T}{(\ln x)^{\frac{1}{2}} |\zeta(\frac{1}{2}+it)|^{\frac{1}{2}} \sqrt{x} \frac{dt}{t+1}}\right| \ll
		(\ln x)^{\frac{1}{2}+\frac{5}{4}} \sqrt{x} = (\ln x)^{\frac{7}{4}} \sqrt{x}.$$
		The same estimation holds for $j_7$.
		
		The main term arises from the calculation of $j_4$ and $j_5$.
		As $L(1, \chi_4) = \frac{\pi}{4}$, then, using the proof of the theorem \ref{theorem1}, we get:
		$$F(s) = \frac{e^{-\frac{\pi i}{4}}}{\sqrt[4]{u}} \Pi(u),$$
		where
		$$\Pi(u) = H(1-u)(L(1-u, \chi_4))^{\frac{1}{4}} (\zeta(2-2u))^{\frac{17}{48}} (L(2-2u, \chi_4))^{-\frac{19}{48}} \sqrt[4]{w(1-u)}.$$
		Suppose that $N \ge 0$ is fixed. Then:
		$$\Pi(u)= \Pi_0+\Pi_1u + \Pi_2u^2 + \ldots + \Pi_N u^N + O_N(u^{N+1}).$$
		Then:
		$$j_4 = \frac{1}{2 \pi i} \int_{\frac{1}{2} + \frac{1}{\ln x} + i \cdot 0}^{1+i \cdot 0}{F(\sigma + i \cdot 0) \frac{(x+h)^s- x^s}{s} ds} = $$
		$$= \frac{1}{2\pi i}\int_{x}^{x+h}{\int_{\frac{1}{2}+\frac{1}{\ln x}}^{1}{F(\sigma+i \cdot 0)y^{\sigma-1}d\sigma}\;dy}=
		\frac{e^{-\frac{\pi i}{4}}}{2\pi i} \int_{x}^{x+h}{{\int_{0}^{\frac{1}{2}-\frac{1}{\ln x}}{\frac{\Pi(u)y^{-u}}{\sqrt[4]{u}}}du}\;dy} =$$
		$$= \frac{e^{-\frac{\pi i}{4}}}{2\pi i} \int_{x}^{x+h}{\left(\sum_{0 \le n \le N}{\Pi_n \int_{0}^{\frac{1}{2}-\frac{1}{\ln x}}{\frac{u^{n}y^{-u}}{\sqrt[4]{u}}du}}+ O\left( J\right)\right)dy},$$
		where
		$$J = \int_{0}^{\frac{1}{2}-\frac{1}{\ln x}}{\frac{u^{N+1}y^{-u}}{\sqrt[4]{u}}du \ll \frac{\Gamma\left(N+\frac{7}{4}\right)}{(\ln y)^{N+\frac{7}{4}}}}.$$
		The contribution of the $n$-th term to the sum is equal to
		$$\int_{0}^{\frac{1}{2}-\frac{1}{\ln x}}{\frac{u^n y^{-u}}{\sqrt[4]{u}}du} =\frac{\Gamma\left( n+\frac{3}{4}\right)}{(\ln y)^{n+\frac{3}{4}}}+ \frac{\theta e n!}{\sqrt{y} \ln y}.$$
		Consequently, 
		$$j_4 = \frac{h e^{-\frac{\pi i}{4}}}{2\pi i} \sum_{0 \le n \le N}{\frac{\Pi_n \Gamma\left( n+\frac{3}{4}\right)}{(\ln x)^{n+\frac{3}{4}}}}+ O\left( \frac{h}{(\ln x)^{N+\frac{7}{4}}}\right) + O\left( \frac{h^2}{x} \frac{1}{(\ln x)^{N+\frac{7}{4}}}\right),$$
		$$j_5 = -\frac{h e^{-\frac{\pi i}{4}}}{2\pi i} \sum_{0 \le n \le N}{\frac{\Pi_n \Gamma\left( n+\frac{3}{4}\right)}{(\ln x)^{n+\frac{3}{4}}}}+ O\left( \frac{h}{(\ln x)^{N+\frac{7}{4}}}\right) + O\left( \frac{h^2}{x} \frac{1}{(\ln x)^{N+\frac{7}{4}}}\right).$$
		
		By Euler's reflection formula for the Gamma-function,
		$$-(j_4+j_5)=
		-\frac{h}{(\ln x)^{\frac{3}{4}}} \left( \sum_{0 \le n \le N}{\frac{(-1)^n \Pi_n}{\Gamma \left( \frac{1}{4}-n\right) (\ln x)^n}} + O\left( \frac{1}{(\ln x)^{N+1}}\right) + O \left( \frac{h}{x(\ln x)^{N+1}}\right)\right).$$
		
		It remains to estimate the sum
		$$\sum_{|\gamma|<T}{j_{\rho}},\quad{\rm where}\quad \rho = \beta +i\gamma,\;\;j_{\rho}=I_1(\rho)+I_2(\rho).$$
		
		\begin{figure}[tbh]
			\begin{center}
				\includegraphics{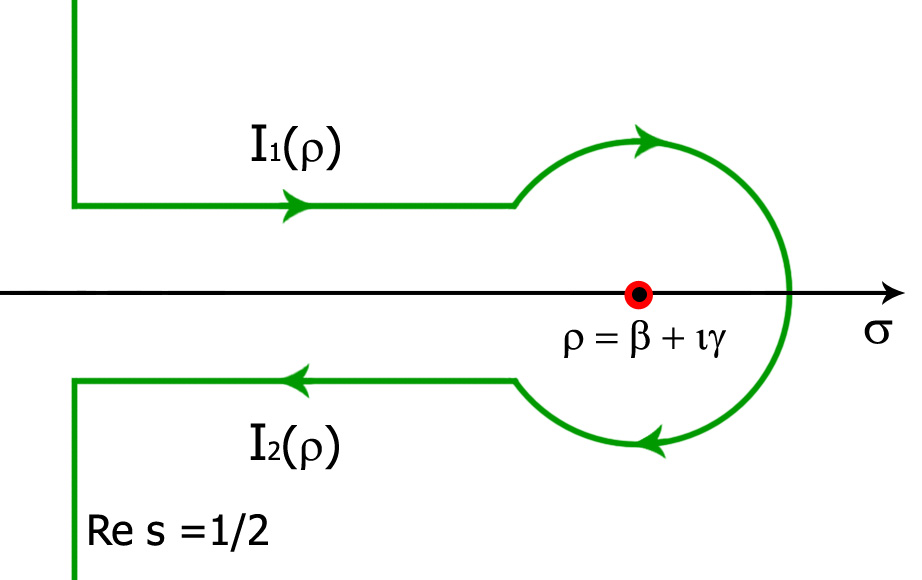}
			\end{center}
		\end{figure}
		
		Firstly, consider the sum over the zeroes of $\zeta(s)$.
		Using the same arguments as above, we get
		$$I_1(\rho) \ll \frac{h}{x} (\ln x)^{\frac{1}{2}} \int_{\frac{1}{2}}^{\beta}{x^{\sigma}  |\zeta(\sigma+i\gamma) L(\sigma+i\gamma, \chi_4)|^{\frac{1}{4}} d\sigma}.$$
		Therefore,
		$$|j_\rho| \ll h {\ln x} \int_{\frac{1}{2}}^{1}{g(\rho,\sigma) \left( \frac{T^{\frac{1}{6}}}{x}\right)^{1-\sigma} d\sigma},$$
		where
		\begin{equation} \label{ggg}
			g(\rho,\sigma) = \left\{
				\begin{aligned}
					1,\;{\rm if}\;\sigma \le \beta,\\ 
					0,\;{\rm if}\;\sigma > \beta.\\
				\end{aligned}
				\right.
		\end{equation}
		Applying the lemma \ref{density},
		$$\sum_{|\gamma|<T}{j_\rho} \ll	 h (\ln x)^{45} \left( \frac{T^{\frac{12}{5} + \frac{1}{6}}}{x}\right)^{1-\varrho (T)}.$$
		
		Now let us evaluate the sum over the zeroes of $L(s, \chi_4)$. We have:
		$$I_1(\rho) \ll h (\ln x)^{\frac{3}{4}+\frac{17}{48}} \int_{\frac{1}{2}}^{\beta}{x^{\sigma-1} T^{\frac{1-\sigma}{6}} d\sigma} \ll
		h (\ln x)^{\frac{53}{48}} \int_{\frac{1}{2}}^{1}{g(\rho; \sigma) \left( \frac{T^{\frac{1}{6}}}{x}\right)^{1-\sigma} d\sigma},$$
		where $g(\rho; \sigma)$ is defined by (\ref{ggg}).
		Applying lemma \ref{density}, we get:
		$$\sum_{|\gamma| < T}{j_{\rho}} \ll h (\ln x)^{\frac{41}{48}} \int_{\frac{1}{2}}^{1-\varrho(T)}{N(\sigma; T, \chi_4) \left( \frac{T^{\frac{1}{12}}}{x}\right)^{1-\sigma}d\sigma} \ll$$
		$$\ll h (\ln x)^{\frac{53}{48}} \int_{\frac{1}{2}}^{1-\varrho(T)}{\sum_{\chi \mmod 4} {N(\sigma; T, \chi_4) \left( \frac{T^{\frac{1}{6}}}{x}\right)^{1-\sigma} d\sigma}} \ll$$
		$$\ll h (\ln x)^{45+\frac{53}{48}} \int_{\frac{1}{2}}^{1-\varrho(T)}{\left( \frac{T^{\frac{1}{6}+\frac{12}{5}}}{x}\right)^{1-\sigma} d\sigma}.$$
		Choosing $T$ from the equation
		$$T^{\frac{12}{5}+\frac{1}{6}} = xD^{-1}(x),\;D(x)=e^{C_1 (\ln x)^{0.8}}$$
		we get
		$$T=x^{\frac{30}{77}} D(x)^{-\frac{30}{77}}.$$
		Now we conclude that the formula
		$$j_0 = \frac{h}{(\ln x)^{\frac{3}{4}}} \sum_{0 \le n \le N}{\frac{(-1)^n \Pi_n}{\Gamma(\frac{1}{4}-n) (\ln x)^n}} + O(J),$$
		$$J \ll \frac{x}{T} (\ln x)^{\frac{17}{4}} + \frac{h}{(\ln x)^{N+\frac{7}{4}}}+ h (\ln x)^{47} \left( \frac{T^{\frac{12}{5} + \frac{1}{6}}}{x}\right)^{1-\varrho (T)},$$
		is asymptotic, if
		$$h = x^{\alpha} e^{C_2 (\ln x)^{0.8}} \gg \frac{x}{T} (\ln x)^2,$$	
		where
		$$\alpha=1-\frac{30}{77}=\frac{47}{77},$$
		which proves the theorem.
	\end{subsection}
\end{section}

\end{document}